\documentclass[twocolumn]{autart}

\usepackage{graphicx}
\usepackage{graphics}
\usepackage{epsfig}
\usepackage{times}
\usepackage{amsmath}
\usepackage{amssymb}
\usepackage{epstopdf}
\usepackage{epsfig}
\usepackage{enumerate}
\usepackage[noadjust]{cite}

\usepackage{makecell}
\usepackage{multirow}
\usepackage{threeparttable}

\newtheorem{Theorem 1}{Theorem}
\newtheorem{Theorem 2}[Theorem 1]{Theorem}
\newtheorem{Theorem 3}[Theorem 1]{Theorem}
\newtheorem{Theorem 4}[Theorem 1]{Theorem}
\newtheorem{Definition 1}{Definition}
\newtheorem{Definition 2}[Definition 1]{Definition}
\newtheorem{Definition 3}[Definition 1]{Definition}
\newtheorem{Definition 4}[Definition 1]{Definition}
\newtheorem{Definition 5}[Definition 1]{Definition}
\newtheorem{Definition 6}[Definition 1]{Definition}
\newtheorem{Definition 7}[Definition 1]{Definition}
\newtheorem{Definition 8}[Definition 1]{Definition}
\newtheorem{Definition 9}[Definition 1]{Definition}
\newtheorem{Remark 1}{Remark}
\newtheorem{Remark 2}[Remark 1]{Remark}
\newtheorem{Remark 3}[Remark 1]{Remark}
\newtheorem{Remark 4}[Remark 1]{Remark}
\newtheorem{Remark 5}[Remark 1]{Remark}
\newtheorem{Lemma 1}{Lemma}
\newtheorem{Lemma 2}[Lemma 1]{Lemma}
\newtheorem{Lemma 3}[Lemma 1]{Lemma}
\newtheorem{Corollary 1}{Corollary}
\newtheorem{Corollary 2}[Corollary 1]{Corollary}
\newtheorem{Corollary 3}[Corollary 1]{Corollary}

\newtheorem{Assumption 1}{Assumption}
\newtheorem{Assumption 2}[Assumption 1]{Assumption}
\newtheorem{Assumption 3}[Assumption 1]{Assumption}
\newtheorem{Proposition 1}{Proposition}
\usepackage{color}

\begin{document}

\begin{frontmatter}
\runtitle{On the Global Synchronization of Pulse-coupled Oscillators Interacting on Chain and Directed Tree Graphs}

\title{On the Global Synchronization of Pulse-coupled Oscillators Interacting on Chain and Directed Tree Graphs\thanksref{footnoteinfo}}

\thanks[footnoteinfo]{The work was supported in part by the National Science Foundation under Grant 1738902.}

\thanks{\bf This paper has been accepted to Automatica as a full paper. Please cite this article as: H. Gao and Y. Wang, On the global synchronization of pulse-coupled oscillators interacting on chain and directed tree graphs. Automatica (2019), https://doi.org/10.1016/j.automatica.2019.02.059.}

\author{Huan Gao}\ead{hgao2@clemson.edu}, 
\author{Yongqiang Wang\thanksref{*}}\ead{yongqiw@clemson.edu}

\thanks[*]{Corresponding author.}

\address{Department of Electrical and Computer Engineering, Clemson University, Clemson, SC 29634, United States}

\begin{keyword}
	Global synchronization; pulse-coupled oscillators; hybrid systems.
\end{keyword}

\begin{abstract}
Driven by increased applications in biological networks and wireless sensor networks, synchronization of pulse-coupled oscillators (PCOs) has gained increased popularity. However, most existing results address the \emph{local} synchronization of PCOs with initial phases constrained in a half cycle, and results on \emph{global} synchronization from \emph{any initial condition} are very sparse. In this paper, we address \emph{global} PCO synchronization from an arbitrary phase distribution under chain or directed tree graphs. Our results differ from existing global synchronization studies on decentralized PCO networks in two key aspects: first, our work allows heterogeneous coupling functions, and we analyze the behavior of oscillators with perturbations on their natural frequencies; secondly, rather than requiring a large enough coupling strength, our results hold under any coupling strength between zero and one, which is crucial because a large coupling strength has been shown to be detrimental to the robustness of PCO synchronization to disturbances.
\end{abstract}

\end{frontmatter}

\section{Introduction}
Pulse-coupled oscillators (PCOs) are limit cycle oscillators coupled through exchanging pulses at discrete time instants. They were originally proposed to model the synchronization phenomena in biological systems, such as contracting cardiac cells, flashing fireflies, and firing neurons \cite{peskin:75, mirollo:90, Ermentrout:96}. Due to their amazing scalability, simplicity, and robustness, recently they have found applications in wireless sensor networks \cite{hong:05, Pagliari:11, Hu:06, Simeone:08}, image processing \cite{rhouma2001self}, and motion coordination \cite{Huan_Tac_2017}.

Early results on PCO synchronization were motivated by biological applications, and normally assume a fixed interaction or coupling mechanism \cite{peskin:75, mirollo:90}. In engineering applications, such restrictions do not exist any more. In fact, the interaction mechanism becomes a design variable that provides opportunities to achieve desired performance. For example, \cite{mauroy2011dichotomic} and \cite{nishimura2011robust} designed the interaction to improve the robustness to communication delays. Our prior work \cite{wang_tsp2:12} optimized the interaction, i.e., phase response function (PRF), to improve the speed of synchronization. However, most of these results are for local synchronization assuming that the initial phases are restricted within a half cycle \cite{Vreeswijk:94, Ernst:95, hansel1995synchrony, Ermentrout:96, Kirk:97, acker2003synchronization, achuthan2009phase, canavier2010pulse, nishimura2011robust, lamar2010effect, timme2002coexistence, timme2008simplest, memmesheimer2010stable, kannapan2016synchronization, wang_tsp2:12, wang_tsp:12, wang_tsp:13, proskurnikov2015event, goel2002synchrony, proskurnikov2016synchronization, dror1999mathematical}.

\begin{table*}
	\centering{Table 1. Comparison our results with other results.}
	\begin{center}
		\begin{threeparttable}	
			\renewcommand\arraystretch{2}
			\begin{tabular}{|c |c |c |c |c |c |}
				\hline
				\multicolumn{2}{|c|}{\multirow{2}{*}{\makecell{ }}} & \multicolumn{2}{c|}{\makecell{ Homogeneous coupling \\ }} & \multicolumn{2}{c|}{\makecell{Heterogeneous coupling \\ }} \\ \cline{3-6}
				\multicolumn{2}{|c|}{} & \makecell{PCO network having \\ (at least) a global node\tnote{1}} & \makecell{Decentralized \\ PCO networks} & \makecell{PCO network having \\ (at least) a global node} & \makecell{Decentralized \\ PCO networks} \\ \hline
				\multirow{3}{*}{\makecell{\\ \\ Non-global \\ synchronization}} & \makecell{Local \\ synchronization} & \makecell{\cite{Vreeswijk:94, Ernst:95, hansel1995synchrony, Ermentrout:96, Kirk:97, goel2002synchrony, acker2003synchronization, konishi:08, achuthan2009phase, canavier2010pulse, nunez2015synchronization}} & \makecell{\cite{goel2002synchrony, nishimura2011robust, lamar2010effect, timme2002coexistence, timme2008simplest, memmesheimer2010stable, proskurnikov2015event, kannapan2016synchronization, wang_tsp2:12, wang_tsp:12, wang_tsp:13, nunez2015synchronization, nunez2015global}} & \makecell{\cite{nunez2016synchronization}} & \makecell{\cite{dror1999mathematical, proskurnikov2016synchronization}} \\ \cline{2-6}	
				& \makecell{Almost global \\ synchronization \\ or synchronization \\ with probability one} & \makecell{\cite{mirollo:90, chen:94, Mathar:96}} & \makecell{\cite{klinglmayr2012guaranteeing, klinglmayr2017convergence, lyu2015synchronization}} & \makecell{$\diagdown$} & \makecell{$\diagdown$} \\ \hline
				\multirow{2}{*}{\makecell{Global \\ synchronization}} & \makecell{Discrete state \\ synchronization} & \cite{an2011nonidentical} & \makecell{\cite{lyu2015synchronization}} & $\diagdown$ & $\diagdown$	\\ \cline{2-6}
				& \makecell{(Continuous) phase \\ synchronization} & \makecell{\cite{konishi:08, klinglmayr2012self, nunez2015synchronization, canavier2017globally}} & \makecell{[32,\,33,\,43\tnote{2} ]} & \makecell{\cite{nunez2016synchronization}} & \makecell{\bf This paper} \\ \hline				
			\end{tabular}
		\end{threeparttable}
	\end{center}
	\begin{tablenotes}
		\item [] $^{1}$ A node is called as a global node if it is directly connected to all the other nodes.	
		\item [] $^{2}$ Note that when the maximum degree of an undirected tree graph is not over $3$, \cite{lyu2018global} obtained global synchronization results for the conventional phase-only PCO model, though results were also obtained under general undirected tree graphs for a more complicated PCO model with multiple additional state variables.
	\end{tablenotes}
\end{table*}

Assuming restricted initial phase distribution severely hinders the application of PCO based synchronization, since in distributed systems it is hard to control the initial phase distribution. Recently, efforts have emerged to address global PCO synchronization from an arbitrary initial phase distribution. However, these results focus on special graphs, such as all-to-all graph \cite{konishi:08, klinglmayr2012self, canavier2017globally, nunez2015synchronization}, cycle graph \cite{nunez2015global}, strongly-rooted graph \cite{nunez2015synchronization}, or master/slave graph \cite{nunez2016synchronization}. Moreover, they rely on sufficiently large coupling strengths, which may not be desirable as large coupling strengths are detrimental to robustness to disturbances \cite{hong:05}.

In this paper, we address the global synchronization of PCOs under arbitrary initial conditions and heterogeneous coupling functions (PRFs). Our main focus is on the global synchronization of PCOs under undirected chain graphs, but the results are easily extendable to PCO synchronization under directed chain/tree graphs. Note that the chain or directed tree graphs are basic elements for constructing more complicated graphs and are desirable in engineering applications where reducing the number of connections is important to save energy consumption and cost in deployment/maintenance. Furthermore, the chain graph has been regarded as the worst-case scenario for synchronization due to its minimum number of connections \cite{klinglmayr2010synchronization}. We also consider oscillators with perturbations on their natural frequencies. Compared with existing results including our prior work (cf. Table 1), this paper has the following contributions: 1) Different from most existing results which focus on local PCO synchronization and assume that the initial phases of oscillators are restricted within a half cycle, our work addresses global synchronization from an arbitrary initial phase distribution; 2) Different from existing global synchronization studies on decentralized PCO networks, our work allows heterogeneous phase response functions, and we analyze the behavior of oscillators with perturbations on their natural frequencies. These scenarios, to our knowledge, have not been considered in any existing global synchronization results on decentralized PCO networks; 3) In contrast to existing global PCO synchronization results requiring a strong enough coupling strength, our results guarantee global synchronization under any coupling strength between zero and one, which is more desirable since a very strong coupling strength, although can bring fast convergence, has been shown to be detrimental to the robustness of synchronization to disturbances \cite{hong:05}.

It is worth noting that even in the theoretical derivation point of view, this paper also differs significantly from our prior work \cite{ nunez2015synchronization, nunez2015global, nunez2016synchronization}: 1) Different from our prior work \cite{ nunez2015synchronization, nunez2015global, nunez2016synchronization} whose proofs are essentially based on local synchronization analysis, this work presents a direct global analyzing approach. More specifically, to obtain global synchronization results, our prior work \cite{ nunez2015synchronization, nunez2015global, nunez2016synchronization} used {\bf strong enough} coupling strengths to reduce the network to a state where all phases are contained in a half cycle, and then achieved global synchronization based on local synchronization analysis. In comparison, this work studies the systematic evolution of phases even when they are not restricted in a half cycle, and hence can allow the coupling strength to be any value between zero and one; 2) Although the Lyapunov candidate function seems similar to the one used in our prior work \cite{nunez2015global}, the analysis here is much more complicated due to the considered more complicated scenarios (arbitrary coupling strength between zero and one and heterogeneous PRFs). In fact, to address synchronization under such scenarios, we had to introduce Invariance Principle, which is not needed in our prior results \cite{ nunez2015synchronization, nunez2015global, nunez2016synchronization} due to their simple dynamics brought by strong and homogeneous coupling.

The outline of this paper is as follows. Section 2 introduces preliminary concepts. A hybrid model for PCO networks and its dynamical properties are presented in Section 3. In Section 4, we analyze global synchronization on both chain and directed tree graphs and provide robustness analysis under frequency perturbations. Numerical experiments are given in Section 5. Finally, we conclude the paper in Section 6.


\section{Preliminaries}

\subsection{Basic Notation}

$\mathbb{R}$, $\mathbb{R}_{\geq 0}$, and $\mathbb{Z}_{\geq 0}$ denote real numbers, nonnegative real numbers, and nonnegative integers, respectively. $\mathbb{R}^n$ denotes the Euclidean space of dimension $n$, and $\mathbb{R}^{n \times n}$ denotes the set of $n \times n$ square matrices with real coefficients. $\mathbb{B}$ denotes the closed unit ball in the Euclidean norm. A set-valued map $\it{M} : A \rightrightarrows B$ associates an element $\alpha \in A$ with a set $\it{M}(\alpha) \subseteq B$; the graph of $\it{M}$ is defined as ${\rm{graph}}(\it{M}) := \{(\alpha,\beta)\in A \times B : \beta \in \it{M}(\alpha)\}$. $\it{M}$ is outer-semicontinuous if and only if its graph is closed \cite{variational_book1}. The range of a function $f : \mathbb{R}^n \rightarrow \mathbb{R}^m$ is denoted as ${\rm{rge}}\, f$. The closure of set $\mathcal{A}$ is denoted as $\overline {\mathcal{A}}$. The distance of a vector $x\in \mathbb{R}^n$ to a closed set $\mathcal{A} \subset \mathbb{R}^n$ is denoted as $|x|_\mathcal{A}=\inf_{y\in \mathcal{A}}|x-y|$. The $\mu$-level set of function $V : {\rm{dom}}\, V \rightarrow \mathbb{R}$ is denoted as $V^{-1}(\mu)=\{x\in {\rm{dom}}\, V: V(x)=\mu\}$ \cite{hybrid_book1}.


\subsection{Hybrid Systems}

We use hybrid systems framework with state $ x \in \mathbb{R}^n$ \cite{hybrid_book1}
\begin{equation}\label{hybrid_systems}
\mathcal{H} : \left\lbrace \begin{aligned}
& \dot{x} = f(x), & & x \in \mathcal{C} \\
& x^+ \in G(x), & & x \in \mathcal{D} \\
\end{aligned}\right.
\end{equation}
where $f$, $\mathcal{C}$, $G$, and $\mathcal{D}$ are the {\it flow map}, {\it flow set}, {\it jump map}, and {\it jump set}, respectively. The hybrid system can be represented by $\mathcal{H}=(\mathcal{C},f,\mathcal{D},G)$. In hybrid system, a {\it hybrid time point} $(t, \, j)\in E$ is parameterized by both $t$, the amount of time passed since initiation, and $j$, the number of jumps that have occurred. A subset $E \subset \mathbb{R}_{\geq 0} \times \mathbb{Z}_{\geq 0}$ is a {\it hybrid time domain} if it is the union of a finite or infinite sequence of interval $[t_k,t_{k+1}] \times \{k\}$. A {\it solution} to $\mathcal{H}$ is a function $\phi: E \rightarrow \mathbb{R}^n$ where $\phi$ satisfies the dynamics of $\mathcal{H}$, $E$ is a hybrid time domain, and for each $j \in \mathbb{N}$, the function $t \mapsto \phi(t,j)$ is locally absolutely continuous on $I_j=\{t:(t,j)\in E\}$. $\phi(t,j)$ is called a {\it hybrid arc}. A hybrid arc $\phi$ is {\it nontrivial} if its domain contains at least two points, is {\it maximal} if it is not the truncation of another solution, and is {\it complete} if its domain is unbounded. Moreover, a hybrid arc $\phi$ is {\it Zeno} if it is complete and $\sup_t \rm{dom} \, \phi < \infty$, is {\it continuous} if it is nontrivial and $\rm{dom} \, \phi \subset \mathbb{R}_{\geq 0} \times \{0\}$, is {\it eventually continuous} if $J=\sup_j \rm{dom} \, \phi < \infty$ and $\rm{dom} \, \phi \cap (\mathbb{R}_{\geq 0} \times \{J\})$ contains at least two points, is {\it discrete} if it is nontrivial and $\rm{dom} \, \phi \subset \{0\} \times \mathbb{Z}_{\geq 0}$, and is {\it eventually discrete} if $T=\sup_t \rm{dom} \, \phi < \infty$ and $\rm{dom} \, \phi \cap (\{T\} \times \mathbb{Z}_{\geq 0})$ contains at least two points. Given a set $\mathcal{M}$, we denote $\mathcal{S_H(M)}$ the set of all maximal solutions $\phi$ to $\mathcal{H}$ with $\phi(0, 0) \in \mathcal{M}$.

Some notions and results for the hybrid system $\mathcal{H}$ from \cite{hybrid_book1} which will be used in this paper are given as follows.

\begin{Definition 1}\label{hybrid_basic}
	$\mathcal{H}=(\mathcal{C},f,\mathcal{D},G)$ satisfies the hybrid basic conditions if: 1) $\mathcal{C}$ and $\mathcal{D}$ are closed in $\mathbb{R}^n$; 2) $f : \mathbb{R}^n \rightarrow \mathbb{R}^n $ is continuous and locally bounded on $\mathcal{C} \subset {\rm{dom}} \, f$; and 3) $G : \mathbb{R}^n \rightrightarrows \mathbb{R}^n$ is outer-semicontinuous and locally bounded on $\mathcal{D} \subset {\rm{dom}} \, G$.
\end{Definition 1}

\begin{Definition 1}\label{strongly forward invariant}
	A set $S\subset \mathbb{R}^n$ is said to be strongly forward invariant if for every $\phi \in \mathcal{S_H}(S)$, $\rm{rge} \, \phi \subset S$.
\end{Definition 1}

\begin{Definition 3}\label{pre-forward completeness}
	Given a set $S\subset \mathbb{R}^n$, a hybrid system $\mathcal{H}$ on $\mathbb{R}^n$ is pre-forward complete from $S$ if every $\phi \in \mathcal{S_H}(S)$ is either bounded or complete.
\end{Definition 3}

\begin{Definition 1}\label{uniformly attractive}
	A compact set $\mathcal{A} \subset \mathbb{R}^n$ is said to be uniformly attractive from a set $S \subset \mathbb{R}^n$ if every $\phi \in \mathcal{S_H}(S)$ is bounded and for every $\varepsilon >0$ there exists $\tau>0$ such that $|\phi (t, \, j)|_{\mathcal{A}} \leq \varepsilon$ for every $\phi \in \mathcal{S_H}(S)$ and $(t, \, j) \in \rm{dom} \, \phi$ with $t+j \geq \tau$.
\end{Definition 1}

\begin{Definition 2}\label{locally asymptotically stable}
	A compact set $\mathcal{A} \subset \mathbb{R}^n$ is said to be
	\begin{itemize}
		\item stable for $\mathcal{H}$ if for every $\varepsilon>0$ there exists $\delta>0$ such that every solution $\phi$ to $\mathcal{H}$ with $|\phi(0,0)|_\mathcal{A}\leq \delta$ satisfies $|\phi(t,j)|_\mathcal{A}\leq \varepsilon$ for all $(t,j)\in {\rm{dom}} \, \phi$;
		\item locally attractive for $\mathcal{H}$ if every maximal solution to $\mathcal{H}$ is bounded and complete, and there exists $\mu>0$ such that every solution $\phi$ to $\mathcal{H}$ with $|\phi(0,0)|_\mathcal{A}\leq \mu$ converges to $\mathcal{A}$, i.e., $\lim_{t+j\rightarrow \infty}|\phi(t,j)|_\mathcal{A}=0$ holds;
		\item locally asymptotically stable for $\mathcal{H}$ if it is both stable and locally attractive for $\mathcal{H}$.
	\end{itemize}
\end{Definition 2}

\begin{Definition 3}\label{basin of attraction}
	Let $\mathcal{A} \subset \mathbb{R}^n$ be locally asymptotically stable for $\mathcal{H}$. Then the basin of attraction of $\mathcal{A}$, denoted by $\mathcal{B_A}$, is the set of points such that every $\phi \in \mathcal{S_H(B_A)}$ is bounded, complete, and $\lim_{t+j\rightarrow \infty}|\phi(t,j)|_\mathcal{A}=0$.
\end{Definition 3}

\begin{Definition 4}\label{tau-epsilon}
	Given $\tau,\varepsilon>0$, two hybrid arcs $\phi_1$ and $\phi_2$ are $(\tau,\varepsilon)$-close if
	\begin{itemize}
		\item $\forall \, (t,j)\in {\rm{dom}} \, \phi_1$ with $t+j\leq \tau$ there exists $s$ such that $(s,j)\in {\rm{dom}} \, \phi_2$, $|t-s|<\varepsilon$ and $|\phi_1(t,j)-\phi_2(s,j)|<\varepsilon$;
		\item $\forall \, (t,j)\in {\rm{dom}} \, \phi_2$ with $t+j\leq \tau$ there exists $s$ such that $(s,j)\in {\rm{dom}} \, \phi_1$, $|t-s|<\varepsilon$ and $|\phi_2(t,j)-\phi_1(s,j)|<\varepsilon$.
	\end{itemize}	
\end{Definition 4}

\begin{Lemma 1}\label{Theorem 8.2}
	(Theorem 8.2 in \cite{hybrid_book1}) Consider a continuous function $V: \mathbb{R}^n \rightarrow \mathbb{R}$, any functions $u_C, u_D: \mathbb{R}^n \rightarrow [-\infty, \, \infty]$, and a set $U \subset \mathbb{R}^n$ such that $u_C(z) \leq 0$, $u_D(z) \leq 0$ for every $z \in U$ and such that the growth of $V$ along solutions to $\mathcal{H}$ is bounded by $u_C$, $u_D$ on $U$. Let a precompact solution $\phi^* \in \mathcal{S_H}$ be such that $\overline{\rm{rge} \, \phi^*} \subset U$. Then, for some $r \in V(U)$, $\phi^*$ approaches the nonempty set that is the largest weakly invariant subset of $V^{-1}(r) \cap U \cap \big[\overline{u_C^{-1}(0)} \cup \big(u_D^{-1}(0) \cap G(u_D^{-1}(0))\big)\big]$.
\end{Lemma 1}
\begin{Lemma 1}\label{Proposition 7.5}
	(Proposition 7.5 in \cite{hybrid_book1}) Let $\mathcal{H}$ be nominally well-posed. Suppose that a compact set $\mathcal{A} \subset \mathbb{R}^n$ has the following properties: 1) it is strongly forward invariant, and 2) it is uniformly attractive from a neighborhood of itself, i.e., there exists $\mu >0$ such that $\mathcal{A}$ is uniformly attractive from $\mathcal{A}+ \mu \mathbb{B}$. Then the compact set $\mathcal{A}$ is locally asymptotically stable.
\end{Lemma 1}
\begin{Lemma 1}\label{Proposition 6.34}
	(Proposition 6.34 in \cite{hybrid_book1}) Let $\mathcal{H}$ be well-posed. Suppose that $\mathcal{H}$ is pre-forward complete from a compact set $K \subset \mathbb{R}^n$ and $\rho: \mathbb{R}^n \rightarrow \mathbb{R}_{\geq 0}$. Then for every $\varepsilon >0$ and $\tau \geq 0$, there exists $\delta >0$ with the following property: for every solution $\phi_\delta$ to $\mathcal{H}_{\delta\rho}$ with $\phi_\delta (0, \, 0) \in K+\delta \mathbb{B}$, there exists a solution $\phi$ to $\mathcal{H}$ with $\phi(0, \, 0) \in K$ such that $\phi_\delta$ and $\phi$ are $(\tau,\varepsilon)$-close.
\end{Lemma 1}


\subsection{Communication Graph}

We use a graph $\mathcal{G}=(\mathcal{V}, \, \mathcal{E}, \, \mathcal{W})$ to represent the interaction pattern of PCOs, where the node set $\mathcal{V} =\{1,2,\ldots,N\}$ denotes all oscillators. $\mathcal{E} \subseteq \mathcal{V} \times \mathcal{V}$ is the edge set, whose elements are such that $(i, \, j) \in \mathcal{E}$ holds if and only if node $j$ can receive messages from node $i$. We assume that no self edge exists, i.e., $(i, \, i) \notin \mathcal{E}$. $\mathcal{W} =[w_{ij}]\in \mathbb{R}^{N \times N}$ is the weighted adjacency matrix of $\mathcal{G}$ with $w_{ij} \geq 0$, where $w_{ij}>0$ if and only if $(i, \, j) \in \mathcal{E}$ holds. The out-neighbor set of node $i$, which represents the set of nodes that can receive messages from node $i$, is denoted as $\mathcal{N}_i^{out} := \{j \in \mathcal{V} : (i, \, j) \in \mathcal{E}\}$.

We focus on chain graphs (both undirected and directed) and directed tree graphs which are defined as follows:
\begin{Definition 5}
	An undirected chain graph $\mathcal{G}$ is a graph whose nodes can be indexed such that there exist two edges $(i, \, i+1)$ and $(i+1, \, i)$ between nodes $i$ and $i+1$ for $i=1,2,\ldots,N-1$.
\end{Definition 5}
\begin{Definition 6}
	A directed chain graph $\mathcal{G}$ is a graph whose nodes can be indexed such that there is only one edge between nodes $i$ and $i+1$ for $i=1,2,\ldots,N-1$ and all edges are directed in the same direction. Without loss of generality, we suppose that the edge between nodes $i$ and $i+1$ is $(i, \, i+1)$.
\end{Definition 6}
\begin{Definition 7}\label{tree_graph}
	A directed tree graph $\mathcal{G}$ is a cycle-free graph with a designated node as a root such that the root has exactly one directed chain to every other node.
\end{Definition 7}

\begin{figure}
	\begin{center}
	\includegraphics[width=0.31\textwidth]{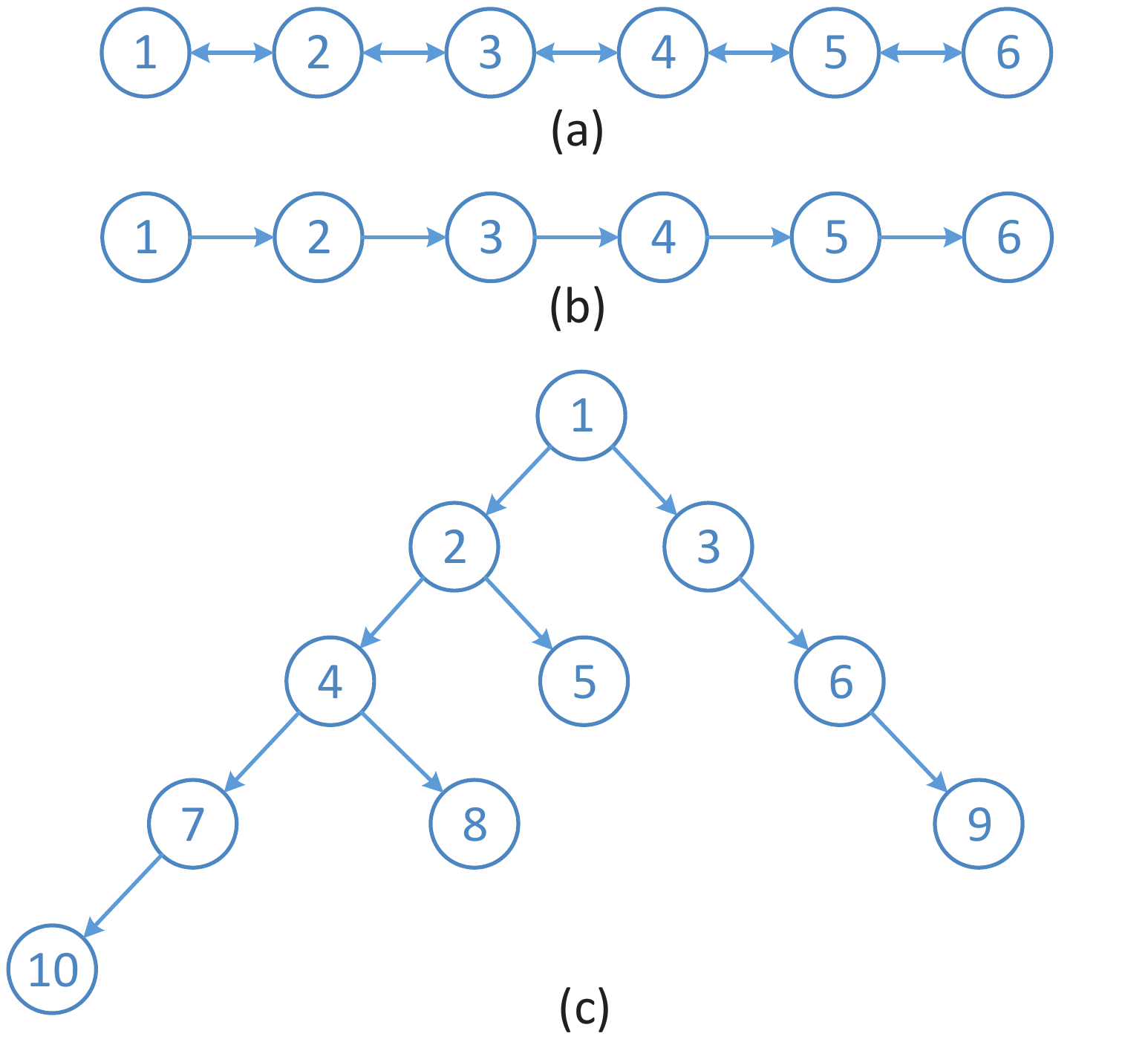}
	\end{center}\caption{Illustration of graphs: (a) undirected chain graph with six nodes; (b) directed chain graph with six nodes; (c) directed tree graph with ten nodes.}
	\label{fig:1}
\end{figure} 

Examples of undirected chain graph, directed chain graph, and directed tree graph are given in Fig. \ref{fig:1}.

\section{Problem Statement}

\subsection{System Description}

We consider $N$ PCOs interacting on a graph $\mathcal{G}=(\mathcal{V}, \, \mathcal{E}, \, \mathcal{W})$. Each oscillator is characterized by a phase variable $x_i \in [0, \, 2\pi]$ for each $i \in \mathcal{V}$. Each phase variable $x_i$ evolves from $0$ to $2\pi$ according to integrate-and-fire dynamics, i.e., $\dot{x_i}=\omega$, where $\omega \in \mathbb{R}_{> 0}$ is the natural frequency of the oscillators. When $x_i$ reaches $2\pi$, oscillator $i$ fires (emits a pulse) and resets $x_i$ to $0$, after which the cycle repeats. When a neighboring oscillator $j$ receives the pulse from oscillator $i$, it shifts its phase according to its coupling strength $l_j \in (0, \, 1)$ (a scalar value) and its phase response function (PRF) $F_j$ \cite{Ermentrout:96, hansel1995synchrony, dror1999mathematical, Izhikevich:07, achuthan2009phase, canavier2010pulse}, i.e., $x_j^+ = x_j+l_jF_j(x_j)$, where $x_j^+$ denotes the phase right after phase shift.


\subsection{Hybrid Model and Dynamical Properties}

Due to the hybrid behavior of PCOs similar to \cite{nunez2015synchronization, nunez2015global, ferrante2017robust}, we model them as a hybrid system $\mathcal{H}$ with state $x = [x_1,\ldots,x_N]^T$. To this end, we define the flow set $\mathcal{C}$ and the flow map $f(x)$ as follows
\begin{equation}\label{flow}
\mathcal{C}=[0, \, 2\pi]^N, \quad f(x)=\omega \mathbf{1}_N \quad \forall \, x \in \mathcal{C}
\end{equation}
According to \cite{nunez2015synchronization, nunez2015global}, the jump set $\mathcal{D}$ and the jump map $G(x)$ can be defined as the union of the individual jump sets $\mathcal{D}_i$ and individual jump maps $G_i(x)$, respectively
\begin{equation}\label{jump}
\mathcal{D} := \bigcup_{i \in \mathcal{V}} \mathcal{D}_i, \quad G(x) := \bigcup_{i \in \mathcal{V}: x \in \mathcal{D}_i} G_i(x)
\end{equation}
where $\mathcal{D}_i$ is defined as $\mathcal{D}_i = \{x \in \mathcal{C} : x_i=2\pi \}$ and $\forall \, x \in \mathcal{D}_i$, $G_i(x)$ is given by
\begin{equation}\label{individual_jump}
\begin{aligned}
& G_i(x)= \{x^+ : x_i^+=0, \, x_j^+ \in x_j+ w_{ij} F_j(x_j) \ \forall \, j \neq i \}
\end{aligned}
\end{equation}
Note $w_{ij}=l_j\in (0, \, 1)$ if $j \in \mathcal{N}_i^{out}$; otherwise, $w_{ij}=0$.

To make $\mathcal{H}$ an accurate description of PCOs, we make the following assumptions on the PRF $F_j$.
\begin{Assumption 1}
	The graph of $F_j$ for $j \in \mathcal{V}$ is such that ${\rm{graph}}(F_j) \subseteq \{(x_j, \, y_j) : x_j \in [0, \, 2\pi], -x_j \leq y_j \leq 2\pi-x_j \}$.
\end{Assumption 1}

This assumption ensures that $G(\mathcal{D}) \subset \mathcal{C} \cup \mathcal{D} = \mathcal{C}$ since $l_j\in (0, \, 1)$ holds, which avoids the existence of solutions ending in finite time due to jumping outside $\mathcal{C}$. 
\begin{Assumption 2}
	The PRF $F_j$ for $j \in \mathcal{V}$ is an outer-semicontinuous set-valued map with $F_j(0)=F_j(2\pi)=0$.
\end{Assumption 2} 

The constraint $F_j(0)=F_j(2\pi)=0$ rules out discrete and eventually discrete solutions, meaning that PCOs will not fire continuously without rest \cite{nunez2015synchronization, nunez2016synchronization}. In fact, there are at most $N$ consecutive jumps with no flow in between because an incoming pulse cannot trigger an oscillator who just fired to fire again under the constraint $F_j(0)=F_j(2\pi)=0$.

The dynamical properties of $\mathcal{H}$ are characterized as follows.
\begin{Proposition 1}\label{Proposition dynamical properties}
	Under Assumptions 1 and 2, we have
	\begin{enumerate}[1)]
		\item $\mathcal{H}$ satisfies the hybrid basic conditions in Definition \ref{hybrid_basic};
		\item For every initial condition $\xi \in \mathcal{C} \cup \mathcal{D}=\mathcal{C}$, there exists at least one nontrivial solution to $\mathcal{H}$. In particular, every solution $\phi \in \mathcal{S_H(C)}$ is maximal, complete, and non-Zeno;
		\item For every solution $\phi \in \mathcal{S_H(C)}$, $\sup_j {\rm{dom}} \, \phi=\infty$ holds, which rules out the existence of continuous and eventually continuous solutions.
	\end{enumerate}
\end{Proposition 1}

{\it Proof}: First we prove statement 1). According to the hybrid model in (\ref{flow})--(\ref{individual_jump}), $\mathcal{C}$ and $\mathcal{D}$ are closed, and $f$ is continuous and locally bounded on $\mathcal{C}$. Also $G$ is locally bounded since the PRF $F_j$ satisfies Assumption 1. To prove $G$ is outer-semicontinuous on $\mathcal{D}$, it suffices to show that ${\rm{graph}}(G) = \bigcup_{i \in \mathcal{V}} \{(x, \, x^+) : x \in \mathcal{D}_i, \, x^+ \in G_i(x)\}$ is closed. According to \cite{ nunez2015synchronization, nunez2015global, nunez2016synchronization}, the outer-semicontinuity of $F_j$ in Assumption 2 ensures that $\{(x, \, x^+) : x \in \mathcal{D}_i, \, x^+ \in G_i(x)\}$ is closed for $i \in \mathcal{V}$, and hence $G$ is outer-semicontinuous on $\mathcal{D}$. Therefore, $\mathcal{H}$ satisfies the hybrid basic conditions in Definition \ref{hybrid_basic}.

Next we prove statement 2). Since $\mathcal{H}$ satisfies the hybrid basic conditions, according to Proposition $6.10$ in \cite{hybrid_book1}, there exists at least one nontrivial solution to $\mathcal{H}$ for every initial condition $\xi \in \mathcal{C} \cup \mathcal{D}=\mathcal{C}$, and every solution $\phi \in \mathcal{S_H(C)}$ is complete due to the facts that $G(\mathcal{D}) \subset \mathcal{C} \cup \mathcal{D} = \mathcal{C}$ holds and $\mathcal{C}$ is compact, which also implies that $\phi$ is maximal. Since $G(\mathcal{D}) \subset \mathcal{C}$ holds, we have ${\rm{rge}} \, \phi \subset \mathcal{C}$ for every $\phi \in \mathcal{S_H(C)}$. So, according to Definition \ref{strongly forward invariant}, $\mathcal{C}$ is strongly forward invariant. Since the constraint $F_j(0)=F_j(2\pi)=0$ in Assumption 2 rules out complete discrete solutions, from Proposition $6.35$ in \cite{hybrid_book1} we have that $\mathcal{S_H(C)}$ is uniformly non-Zeno, which means that every $\phi \in \mathcal{S_H(C)}$ is non-Zeno.

Finally we prove statement 3). Since every $\phi \in \mathcal{S_H(C)}$ is complete and the length of each flow interval is at most $\frac{2\pi}{\omega}$, we have $\sup_j {\rm{dom}} \, \phi=\infty$. So the existence of continuous and eventually continuous solutions is ruled out. \hfill$\blacksquare$ 
\begin{Remark 1}
	As indicated in \cite{nunez2015global}, such hybrid model $\mathcal{H}$ is able to handle multiple simultaneous pulses, i.e., if an oscillator receives multiple pulses simultaneously, it will respond to these pulses sequentially (in whatever order), but the oscillation behavior is the same as if the components of $x$ jumped simultaneously.
\end{Remark 1}


\subsection{General Delay-Advance PRF}

In this paper, we consider general delay-advance PRFs.
\begin{Assumption 3}
	A delay-advance PRF $F_j$ is such that
	\begin{equation}\label{PRF}
	F_j(x_j)=\left\lbrace \begin{aligned}
	& \, F_j^{(1)}(x_j), & & \text{if} \ x_j \in [0, \, \pi) \\
	& \big\{ F_j^{(1)}(\pi), \, F_j^{(2)}(\pi) \big\}, & & \text{if} \ x_j=\pi \\
	& \, F_j^{(2)}(x_j), & & \text{if} \ x_j \in (\pi, \, 2\pi] \\
	\end{aligned}\right.
	\end{equation}
	where $F_j^{(1)}(x_j)$ and $F_j^{(2)}(x_j)$ are continuous functions on $[0, \, \pi]$ and $[\pi, \, 2\pi]$, respectively, and satisfy
	\begin{equation}\label{PRF_detail}
	\left\lbrace \begin{aligned}
	& F_j^{(1)}(0)=0, \ \ F_j^{(1)}(x_j) \in [-x_j, \, 0) \ \ \text{if} \ \ x_j \in (0, \, \pi]\\
	& F_j^{(2)}(2\pi)=0, \ F_j^{(2)}(x_j) \in (0, \, 2\pi-x_j] \ \ \text{if} \ \ x_j \in [\pi, \, 2\pi)\\
	\end{aligned}\right.
	\end{equation}
\end{Assumption 3}

Similar to \cite{nunez2015synchronization, nunez2015global, nunez2016synchronization}, $F_j$ is an outer-semicontinuous set-valued map. Note that oscillators with phases in $(0, \pi)$ will be delayed after receiving a pulse, meaning that their phases will be pushed closer to zero by each pulse received, whereas oscillators with phases in $(\pi, 2\pi)$ will be advanced, meaning that their phases will be pushed toward $2\pi$ by each pulse. If an oscillator has phase $0$ (or $2\pi$) upon receiving a pulse, its phase is unchanged by the pulse.

\begin{figure}
	\begin{center}
		\includegraphics[width=0.5\textwidth]{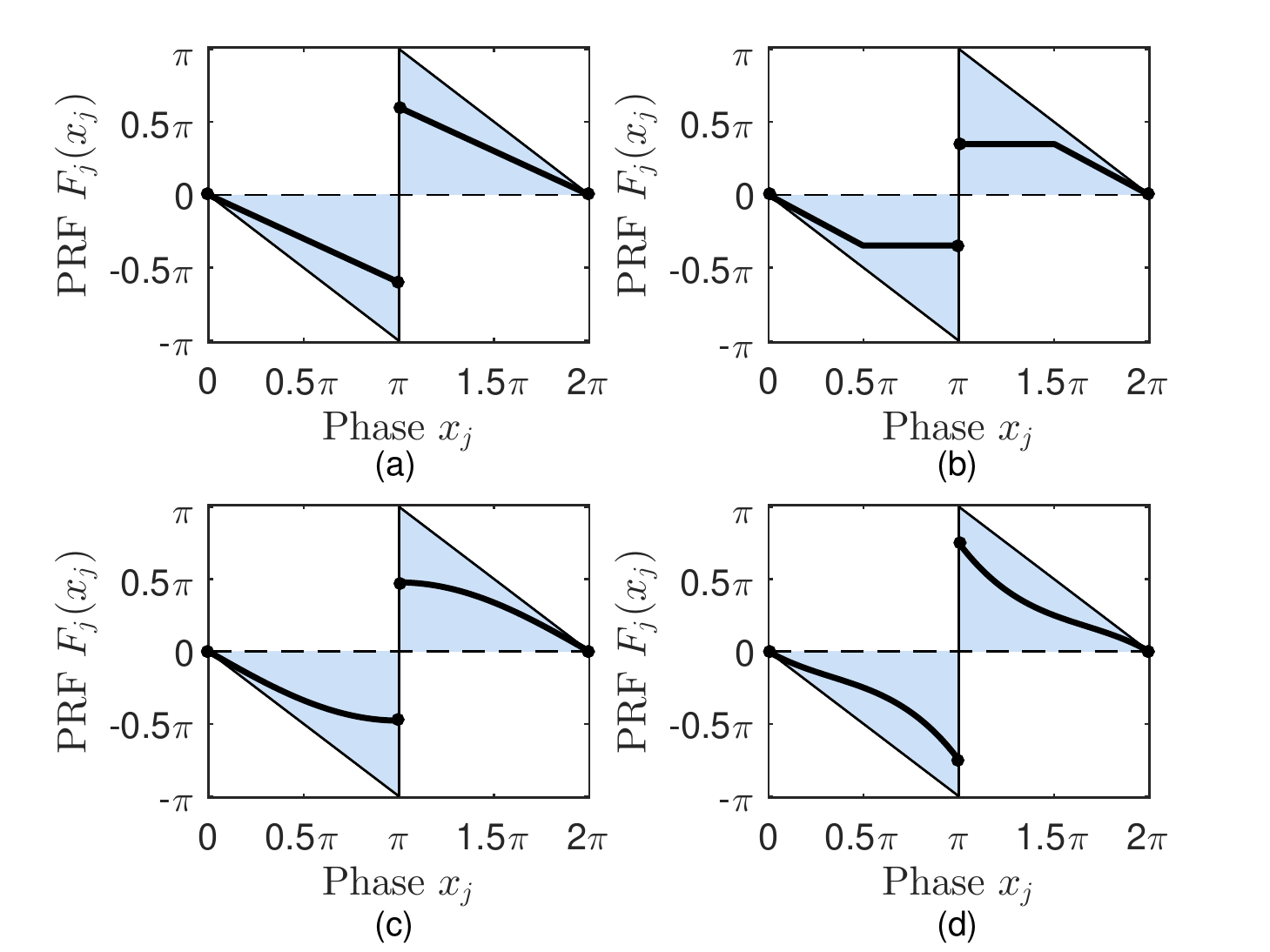}
	\end{center}
	\caption{Examples of the general delay-advance PRF $F_j(x_j)$.}
	\label{fig:2}
\end{figure}

Since Assumption 3 implies Assumptions 1 and 2, the properties of $\mathcal{H}$ in Proposition \ref{Proposition dynamical properties} still hold. Several examples of delay-advance PRF are illustrated in Fig. \ref{fig:2}.

\begin{Remark 1}
	It is worth noting that our PRF can be heterogeneous and is also very general. In fact, it includes the PRFs used in \cite{nunez2015synchronization, nunez2015global, kannapan2016synchronization, wang_tsp2:12, wang_tsp:12, wang_tsp:13, yun2015robustness} as special cases. Therefore, our work has broad potential applications in engineered systems \cite{wang_TCST:12} as well as biological systems \cite{Izhikevich:07}.
\end{Remark 1}

\section{Global Synchronization of PCOs}

In this section, we analyze global PCO synchronization on both chain and directed tree graphs, and provide robustness analysis in the presence of frequency perturbations. 

To this end, we first define the synchronization set $\mathcal{A}$:
\begin{equation}\label{synchronization_set}
\mathcal{A}=\{x \in \mathcal{C}:|x_i-x_j|=0 \, \text{or} \, |x_i-x_j|= 2\pi, \, \forall \, i, \, j \in \mathcal{V}\}
\end{equation}
The PCO network synchronizes if the state $x$ converges to the synchronization set $\mathcal{A}$. Note that $\mathcal{A}$ is compact since it is closed and bounded (included in $\mathcal{C}$ that is bounded).

In the following, we refer to an arc as a connected subset of $[0, \, 2\pi]$ where $0$ and $2\pi$ are associated with each other. So phase difference $\Delta_i$ that measures the length of the shorter arc between $x_i$ and $x_{i+1}$ on the unit cycle is given by
\begin{equation}\label{delta_k}
\Delta_i= \min \{|x_i-x_{i+1}|, \ 2\pi-|x_i-x_{i+1}|\}
\end{equation}
where $x_{N+1}$ is mapped to $x_{1}$ in $\Delta_N$. It is straightforward to show that $\Delta_i$ satisfies $0 \leq \Delta_i \leq \pi$. 

To measure the degree of synchronization, we define $L$ as
\begin{equation}\label{Lyapunov}
L=\sum_{i=1}^{N}\Delta_i
\end{equation}
Since $0 \leq \Delta_i \leq \pi$ holds, we have $0 \leq L \leq N\pi$. Note that both $\Delta_i$ for $i\in \mathcal{V}$ and $L$ are dependent on $x$, and $L$ is positive definite with respect to $\mathcal{A}$ on $\mathcal{C}\cup \mathcal{D}=\mathcal{C}$ because $L=0$ holds if and only if $\Delta_1=\Delta_2=\cdots=\Delta_N=0$ holds. Therefore, in order to prove synchronization, we only need to show that $L$ will converge to $0$. It is worth noting that $L$ is continuous in $x\in \mathcal{C}$ but not differentiable with respect to it.


\subsection{Global Synchronization on Undirected Chain Graphs}

\begin{Lemma 1}\label{Lemma1}
	For $N$ PCOs interacting on an undirected chain, if the PRF $F_j(x_j)$ satisfies Assumption 3 and $l_j\in (0, \, 1)$ holds for all $j \in \mathcal{V}$, then $L$ in (\ref{Lyapunov}) is nonincreasing along any solution $\phi \in \mathcal{S_H(C)}$.
\end{Lemma 1}

{\it Proof}: Since there is no interaction among oscillators during flows and all oscillators have the same natural frequency, we have that $L$ is constant during flows and its dynamics only depends on jumps. Without loss of generality, we assume that at time $(t_i^*,k_i^*)$, we have $x(t_i^*,k_i^*) \in \mathcal{D}_i$, i.e., $x_i(t_i^*,k_i^*)=2\pi$. (In the following, we omit time index $(t_i^*,k_i^*)$ to simplify the notation.) When oscillator $i$ fires and resets its phase to $x_i^+=0$, an oscillator $j \in \mathcal{N}_i^{out}$ has $x_j^+ \in x_j+l_j F_j(x_j)$ but an oscillator $j\notin \mathcal{N}_i^{out}$ still has $x_j^+ = x_j$. 

For the undirected chain graph, we call oscillator $i-1$ as the left-neighbor of oscillator $i$ for $i=2,3,\ldots,N$, and call oscillator $i+1$ as the right-neighbor of oscillator $i$ for $i=1,2,\ldots,N-1$. Upon the firing of oscillator $i$, if the left-neighbor oscillator $i-1$ exists, it will update its phase and affect $\Delta_{i-2}$ and $\Delta_{i-1}$. Note that for $i=2$, $\Delta_{i-2}$ is mapped to $\Delta_N$. Similarly, if the right-neighbor oscillator $i+1$ exists, $\Delta_{i}$ and $\Delta_{i+1}$ will be affected. No other $\Delta_k$s will be affected by this pulse, i.e., $\Delta_{k}^+=\Delta_{k}$ holds for $k \notin \{i-2,i-1,i,i+1\}$ where $\Delta_{k}^+$ denotes the phase difference between oscillators $k$ and $k+1$ after the jump. Therefore, we only need to consider two situations when oscillator $i$ fires, i.e., how $\Delta_{i-2}$ and $\Delta_{i-1}$ change if the left-neighbor oscillator $i-1$ exists and how $\Delta_i$ and $\Delta_{i+1}$ change if the right-neighbor oscillator $i+1$ exists.

\textbf{Situation I:} If the left-neighbor oscillator $i-1$ exists, from (\ref{individual_jump}) and (\ref{PRF}) we have
\begin{equation}\label{phase_update1}
\begin{aligned}
& x_{i-1}^+= \left\lbrace \begin{aligned}
& x_{i-1}+l_{i-1} F_{i-1}^{(1)}(x_{i-1}), & & \text{if} \ x_{i-1}\in [0, \, \pi]\\
& x_{i-1}+l_{i-1} F_{i-1}^{(2)}(x_{i-1}), & & \text{if} \ x_{i-1}\in [\pi, \, 2\pi]\\
\end{aligned}\right.\\
\end{aligned}
\end{equation}
To facilitate the proof, we use an nonnegative variable $\delta_{i-1}$ to denote the jump magnitude of oscillator $i-1$. According to (\ref{PRF_detail}) and $l_{i-1} \in (0, \, 1)$, $\delta_{i-1}$ is determined by
\begin{equation}\label{ddelta_i-1}
\begin{aligned}
\delta_{i-1} = \left\lbrace \begin{aligned}
& -l_{i-1} F_{i-1}^{(1)}(x_{i-1}), & & \text{if} \ x_{i-1}\in [0, \, \pi]\\
& l_{i-1} F_{i-1}^{(2)}(x_{i-1}), & & \text{if} \ x_{i-1}\in [\pi, \, 2\pi]\\
\end{aligned}\right.\\
\end{aligned}
\end{equation}

\begin{figure}
	\begin{center}
		\includegraphics[width=0.49\textwidth]{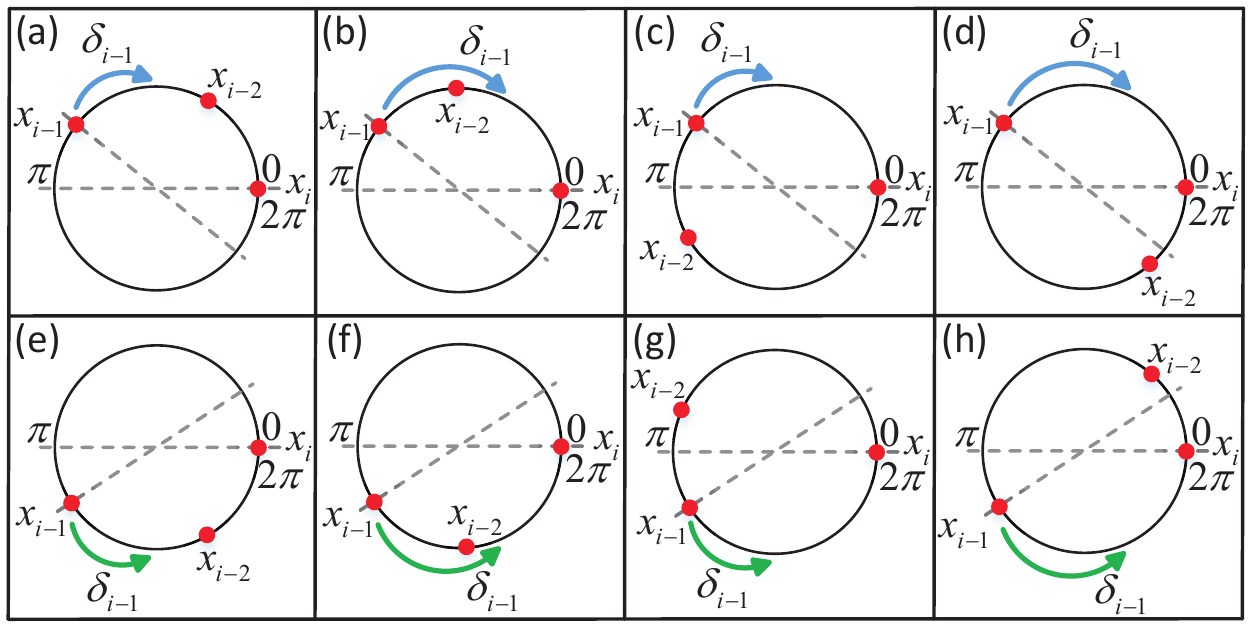}
	\end{center}
	\caption{Illustration of Situation I.}
	\label{fig:Lemma1_fig_1}
\end{figure}

Since $x_i=2\pi$ and $x_i^+=0$ hold, from (\ref{phase_update1}) and (\ref{ddelta_i-1}) we know that oscillator $i-1$ jumps $\delta_{i-1}$ towards oscillator $i$, as illustrated in Fig. \ref{fig:Lemma1_fig_1}. So we have $\Delta_{i-1}^+ = \Delta_{i-1}-\delta_{i-1}$.

Now we analyze how $\Delta_{i-2}$ changes upon oscillator $i$'s firing. Note that $x_{i-2}^+=x_{i-2}$ holds as $i-2 \notin \mathcal{N}_i^{out}$. According to the direction of oscillator $i-1$'s jump and the relationship between $\delta_{i-1}$ and $\Delta_{i-2}$, we have four following cases:

\textbf{Case 1:} If oscillator $i-1$ jumps $\delta_{i-1}$ towards oscillator $i-2$ and $\delta_{i-1} \leq \Delta_{i-2}$ holds (cf. Fig. \ref{fig:Lemma1_fig_1} (a) and (e)), we have $\Delta_{i-2}^+ = \Delta_{i-2}-\delta_{i-1}$, which leads to
\begin{equation}\label{relation_case1}
\Delta_{i-1}^+ + \Delta_{i-2}^+ =\Delta_{i-1} + \Delta_{i-2} -2 \delta_{i-1} \leq \Delta_{i-1} + \Delta_{i-2}
\end{equation}
Note that the equality holds if and only if $\delta_{i-1}=0$ exists, i.e., $\Delta_{i-1}^+ + \Delta_{i-2}^+ =\Delta_{i-1} + \Delta_{i-2}$ holds if and only if $\Delta_{i-2}^+ = \Delta_{i-2}-\delta_{i-1}= \Delta_{i-2}+\delta_{i-1}$ holds.

\textbf{Case 2:} If oscillator $i-1$ jumps $\delta_{i-1}$ towards oscillator $i-2$ and $\delta_{i-1} > \Delta_{i-2}$ holds (cf. Fig. \ref{fig:Lemma1_fig_1} (b) and (f)), we have $\Delta_{i-2}^+ = \delta_{i-1} - \Delta_{i-2}$. So it follows 
\begin{equation}\label{relation_case2}
\begin{aligned}
\Delta_{i-1}^+ + \Delta_{i-2}^+ =\Delta_{i-1} - \Delta_{i-2} \leq \Delta_{i-1} + \Delta_{i-2}
\end{aligned}
\end{equation}
where the equality occurs when $\Delta_{i-2}=0$, i.e., $\Delta_{i-1}^+ + \Delta_{i-2}^+ =\Delta_{i-1} + \Delta_{i-2}$ holds if and only if $\Delta_{i-2}^+ = \delta_{i-1} - \Delta_{i-2}= \Delta_{i-2}+\delta_{i-1}$ holds.

\textbf{Case 3:} If oscillator $i-1$ jumps $\delta_{i-1}$ away from oscillator $i-2$ and $ \Delta_{i-2} + \delta_{i-1} \leq \pi$ holds (cf. Fig. \ref{fig:Lemma1_fig_1} (c) and (g)), we have $\Delta_{i-2}^+ = \Delta_{i-2}+\delta_{i-1}$, which leads to
\begin{equation}\label{relation_case3}
\begin{aligned}
\Delta_{i-1}^+ + \Delta_{i-2}^+ =\Delta_{i-1} + \Delta_{i-2}
\end{aligned}
\end{equation}

\textbf{Case 4:} If oscillator $i-1$ jumps $\delta_{i-1}$ away from oscillator $i-2$ and $ \Delta_{i-2} + \delta_{i-1} > \pi$ holds (cf. Fig. \ref{fig:Lemma1_fig_1} (d) and (h)), we have $\Delta_{i-2}^+ = 2\pi-\Delta_{i-2} -\delta_{i-1} <\pi< \Delta_{i-2}+\delta_{i-1}$ and 
\begin{equation}\label{relation_case4}
\begin{aligned}
\Delta_{i-1}^+ + \Delta_{i-2}^+ & < (\Delta_{i-1}-\delta_{i-1}) + (\Delta_{i-2} + \delta_{i-1}) \\
&= \Delta_{i-1} + \Delta_{i-2}
\end{aligned}
\end{equation}

Summarizing the above four cases, we have
\begin{equation}\label{relation_situation1}
\begin{aligned}
\Delta_{i-1}^+ + \Delta_{i-2}^+ \leq \Delta_{i-1} + \Delta_{i-2}
\end{aligned}
\end{equation}
where the equality occurs when $\Delta_{i-2}^+ = \Delta_{i-2}+\delta_{i-1}$.

\textbf{Situation II:} If the right-neighbor oscillator $i+1$ exists, it will update its phase according to (\ref{individual_jump}) and (\ref{PRF}) as follows
\begin{equation}\label{phase_update_i+1}
\begin{aligned}
& x_{i+1}^+= \left\lbrace \begin{aligned}
& x_{i+1}+l_{i+1} F_{i+1}^{(1)}(x_{i+1}), & & \text{if} \ x_{i+1}\in [0, \, \pi]\\
& x_{i+1}+l_{i+1} F_{i+1}^{(2)}(x_{i+1}), & & \text{if} \ x_{i+1}\in [\pi, \, 2\pi]\\
\end{aligned}\right.\\
\end{aligned}
\end{equation}
Also the nonnegative magnitude of oscillator $i+1$'s phase jump (denoted by $\delta_{i+1}$) is given as
\begin{equation}\label{ddelta_i+1}
\begin{aligned}
\delta_{i+1} = \left\lbrace \begin{aligned}
& -l_{i+1} F_{i+1}^{(1)}(x_{i+1}), & & \text{if} \ x_{i+1}\in [0, \, \pi]\\
& l_{i+1} F_{i+1}^{(2)}(x_{i+1}), & & \text{if} \ x_{i+1}\in [\pi, \, 2\pi]\\
\end{aligned}\right.\\
\end{aligned}
\end{equation}

Since $x_i=2\pi$ and $x_i^+=0$ hold, and oscillator $i+1$ jumps $\delta_{i+1}$ towards oscillator $i$, we have $\Delta_{i}^+ = \Delta_{i}-\delta_{i+1}$.

According to the relationship between $\delta_{i+1}$ and $\Delta_{i+1}$, there are also four cases on the change of $\Delta_{i+1}$. Similar to Situation I, we can obtain the following result
\begin{equation}\label{relation_situation2}
\begin{aligned}
\Delta_{i}^+ + \Delta_{i+1}^+ \leq \Delta_{i} + \Delta_{i+1}
\end{aligned}
\end{equation}
where the equality occurs when $\Delta_{i+1}^+ = \Delta_{i+1}+\delta_{i+1}$.

Summarizing Situation I and Situation II, we can see that $L$ will not increase during jumps. Therefore, $L$ is nonincreasing along any solution $\phi \in \mathcal{S_H(C)}$. \hfill$\blacksquare$ 

Now we are in position to introduce our results for global synchronization on undirected chain graphs.
\begin{thm}
	For $N$ PCOs interacting on an undirected chain, if the PRF $F_j(x_j)$ satisfies Assumption 3 and $l_j\in (0, \, 1)$ holds for all $j \in \mathcal{V}$, then the synchronization set $\mathcal{A}$ in (\ref{synchronization_set}) is globally asymptotically stable, i.e., global synchronization can be achieved from an arbitrary initial condition.
\end{thm}

{\it Proof}: According to the derivation in Lemma \ref{Lemma1}, the continuous function $L$ in (\ref{Lyapunov}) is constant during flows and will not increase during jumps, which implies that $L(g)-L(x)\leq 0$ holds for all $x\in \mathcal{D}$ and $g\in G(x)$. Defining $u_C(x)=0$ for each $x\in \mathcal{C}$ and $u_C(x)=-\infty$ otherwise; $u_D(x)=\max_{g \in G(x)} \{L(g)-L(x)\} \leq 0$ for each $x\in \mathcal{D}$ and $u_D(x)=-\infty$ otherwise, we can bound the growth of $L$ along solutions by $u_C$ and $u_D$ on $\mathcal{C}$ \cite{hybrid_book1}. According to Proposition \ref{Proposition dynamical properties}, every solution $\phi \in \mathcal{S_H(C)}$ is precompact, i.e., complete and bounded, and satisfies $\overline{\rm{rge}\, \phi} \subset \mathcal{C} \cup \mathcal{D} = \mathcal{C}$. From Lemma \ref{Theorem 8.2}, for some $r \in L(\mathcal{C})=[0, \, N\pi]$, $\phi$ approaches the nonempty set that is the largest weakly invariant subset of $L^{-1}(r) \cap \mathcal{C} \cap \big[\overline{u_C^{-1}(0)} \cup \big(u_D^{-1}(0) \cap G(u_D^{-1}(0))\big)\big]$ where $L^{-1}(r)$ denotes the $r$-level set of $L$ defined in Subsection 2.1 (note that Lemma \ref{Theorem 8.2} does not need $L$ to be continuously differentiable in $x\in \mathcal{C}$ \cite{hybrid_book1}). Since $\overline{u_C^{-1}(0)} = \mathcal{C}$ and $u_D^{-1}(0)\cap G(u_D^{-1}(0)) \subset \mathcal{D}$ hold, we have $L^{-1}(r) \cap \mathcal{C} \cap \big[\overline{u_C^{-1}(0)} \cup \big(u_D^{-1}(0) \cap G(u_D^{-1}(0))\big)\big] = L^{-1}(r) \cap \mathcal{C}$.

According to Lemma \ref{Lemma2} in Appendix A, $L$ cannot be retained at any nonzero value along a complete solution $\phi$. So the largest weakly invariant subset of $L^{-1}(r)\cap \mathcal{C}$ is empty for every $r \in (0, \, N\pi]$, which implies that every solution $\phi \in \mathcal{S_H(C)}$ approaches $L^{-1}(0)\cap\mathcal{C}= \mathcal{A}$.

Next we show that $\mathcal{A}$ is locally asymptotically stable. Since every solution $\phi \in \mathcal{S_H(C)}$ approaches $\mathcal{A}$, from Definition \ref{uniformly attractive}, $\mathcal{A}$ is uniformly attractive from $\mathcal{C}$. As Assumption 2 guarantees that $\rm{rge}\, \phi \subset \mathcal{A}$ for every $\phi \in \mathcal{S_H(A)}$, $\mathcal{A}$ is strongly forward invariant according to Definition \ref{strongly forward invariant}. Therefore, from Lemma \ref{Proposition 7.5}, $\mathcal{A}$ is locally asymptotically stable.

To show $\mathcal{A}$ is globally asymptotically stable, it suffices to show that $\mathcal{A}$'s basin of attraction $\mathcal{B_A}$ contains $\mathcal{C}\cup \mathcal{D}=\mathcal{C}$. Since we have shown that the largest weakly invariant subset of $L^{-1}(r)\cap \mathcal{C}$ is empty for every $r \in (0, \, N\pi]$ and every solution $\phi \in \mathcal{S_H(C)}$ approaches $\mathcal{A}$, according to Definition \ref{basin of attraction}, $\mathcal{A}$'s basin of attraction $\mathcal{B_A}$ contains $\mathcal{C}$. Therefore, $\mathcal{A}$ is globally asymptotically stable.

In summary, $\mathcal{A}$ is globally asymptotically stable, meaning that global synchronization can be achieved from an arbitrary initial condition. \hfill$\blacksquare$
\begin{Remark 2}
	Because using four phase differences ($\Delta_{i-2}$, $\Delta_{i-1}$, $\Delta_i$, and $\Delta_{i+1}$, which requires $N \geq 4$) is essential to describe and characterize the dynamics of a general number of $N$ oscillators in a uniform manner, we assumed $N \geq 4$ in the proof. However, the results are also applicable to $N=2$ and $N=3$. In fact, following the analysis in Lemma 4, we can obtain that $L$ is non-increasing when $N=2$ or $3$. Then using the Invariance Principle based derivation in Theorem 1 gives the convergence of $L$ to 0 and thus the achievement of global synchronization for $N=2$ and $3$.
\end{Remark 2}

\begin{Remark 2}
	Compared with existing results in \cite{goel2002synchrony} which show that local synchronization on chain graphs can be obtained as long as the coupling is not too strong, our results can guarantee global synchronization under any coupling strength between zero and one.
\end{Remark 2}

\begin{Remark 3}
	It is worth noting that different from local PCO synchronization analysis \cite{hong:05, goel2002synchrony} and global PCO synchronization analysis under all-to-all topology \cite{nunez2015synchronization, canavier2017globally} where the firing order is time-invariant, the coupling strength $l \in (0, 1)$ cannot guarantee invariant firing order in our considered scenarios, as confirmed by numerical simulations in Fig. \ref{fig:07}.
\end{Remark 3}

\subsection{Global Synchronization on Directed Chain and Tree Graphs}

In this subsection, we extend the global synchronization results to directed chain and tree graphs.

\begin{Corollary 1}
	For $N$ PCOs interacting on a directed chain, if the PRF $F_j(x_j)$ satisfies Assumption 3 and $l_j\in (0, \, 1)$ holds for all $j \in \mathcal{V}$, then the synchronization set $\mathcal{A}$ in (\ref{synchronization_set}) is globally asymptotically stable, i.e., global synchronization can be achieved from an arbitrary initial condition.
\end{Corollary 1}

{\it Proof}: The proof is similar to Theorem 1 and omitted. \hfill$\blacksquare$

\begin{Remark 4}
	Different from the cycle graph in \cite{nunez2015global} where a strong enough coupling strength is required, global synchronization can be achieved here under any coupling strength between zero and one. This is because in the chain case, the absence of interaction between oscillators $1$ and $N$ allows $\Delta_N$ to increase freely until it triggers $L$ to decrease; in other words, the absence of interaction between oscillators $1$ and $N$ breaks the symmetry of the chain graph \cite{golubitsky2003symmetry}, which is key to remove undesired equilibria where $L$ keeps unchanged. In comparison, the symmetry of the cycle graph can make $L$ stay at some undesired equilibria under a weak coupling strength. So a strong enough coupling strength is required in the cycle graph case to achieve global synchronization.
\end{Remark 4}

\begin{Corollary 2}	
	For $N$ PCOs interacting on a directed tree, if the PRF $F_j(x_j)$ satisfies Assumption 3 and $l_j\in (0, \, 1)$ holds for all $j \in \mathcal{V}$, then global synchronization can be achieved from an arbitrary initial condition.
\end{Corollary 2}

{\it Proof}: Suppose in a directed tree graph there are $m$ nodes without any out-neighbors which are represented as $v_1,v_2,\ldots,v_m$. Take the graph in Fig. \ref{fig:1} (c) as an example, nodes $5$, $8$, $9$, and $10$ do not have any out-neighbors. According to Definition \ref{tree_graph}, for every node $v_i$ ($i=1,2,\ldots,m$) there is a unique directed chain from the root $v_r$ to node $v_i$. So the directed tree graph is composed of $m$ directed chains. Note that for every directed chain from the root $v_r$ to node $v_i$, it is not affected by oscillators outside the chain. So the $m$ directed chains are decoupled from each other. According to Corollary 1, global synchronization can be achieved on the directed chain from an arbitrary initial condition if $F_j(x_j)$ satisfies Assumption 3 and if $l_j\in (0, \, 1)$ holds. Adding the fact that the root oscillator $v_r$ belongs to all $m$ directed chains implies synchronization of all PCOs. \hfill$\blacksquare$

\begin{Remark 5}
	Different from the arguments in the proofs of Corollaries 1 and 2, an alternative approach to proving global synchronization on direct chain (and tree) graphs is using inductive reasoning based on the following two facts: first, a parent node can affect its child node but a child node never affects its parent node; secondly, under the given piecewise continuous delay-advance PRF (with values being nonzero in $(0, \, 2\pi)$), the phases of all oscillators on a directed chain will be reduced to within a half cycle, which always leads to synchronization (cf. Theorem 2 in \cite{proskurnikov2016synchronization}).
\end{Remark 5}

\begin{Remark 5}
	Different from the ``probability-one synchronization'' in \cite{klinglmayr2012guaranteeing, klinglmayr2017convergence, lyu2015synchronization} where oscillators synchronize with probability one under a stochastic phase-responding mechanism and the ``almost global synchronization'' in \cite{mirollo:90, chen:94, Mathar:96} where synchronization is guaranteed for all initial conditions except a set of Lebesgue-measure zero, our studied global synchronization is achieved in a deterministic manner from any initial condition, which is not only important theoretically but also mandatory in many safety-critical applications. A typical application justifying the necessity of deterministic global synchronization is synchronization based motion coordination of AUV (autonomous underwater vehicles) \cite{paley2007oscillator} and UAV (unmanned aerial vehicles) \cite{valbuena2014stable}. In such an application, even one single failure in synchronization might be too costly in money, time, energy, or even lives (cf. the multi-UAV based target engagement problem in \cite{furukawa2003coordination}).
\end{Remark 5}

\subsection{Robustness Analysis for Frequency Perturbations}

In this subsection, we analyze the robustness property of PCOs under small frequency perturbations on the natural frequency $\omega$. It is worth noting that robustness is important since frequency perturbations are unavoidable and under an inappropriate synchronization mechanism, even a small difference in natural frequency may accumulate and lead to large phase differences. The hybrid systems model with frequency perturbations is given as follows:
\begin{equation}\label{hybrid_systems_perturbation}
\mathcal{H}_p : \left\lbrace \begin{aligned}
& \dot{x} = \omega \mathbf{1}_N + p, & & x \in \mathcal{C} \\
& x^+ \in G(x), & & x \in \mathcal{D} \\
\end{aligned}\right.
\end{equation}
where $p=[p_1,\ldots,p_N]^T$ represents the frequency perturbations. Using the notion of $(\tau,\varepsilon)$-closeness given in Definition \ref{tau-epsilon} in Subsection 2.2, we have the following result:

\begin{thm}
	Consider $N$ PCOs with frequency perturbations as described by $\mathcal{H}_p$ in (\ref{hybrid_systems_perturbation}). For every $\varepsilon>0$, $\tau \geq 0$, and $\rho : \mathbb{R}^N \rightarrow \mathbb{R}_{\geq 0}$, there exists a scalar $\sigma>0$ such that under any $p \in \sigma \rho(x) \mathbb{B}$ every solution $\phi_p$ to $\mathcal{H}_p$ from $\mathcal{C}$ is $(\tau,\varepsilon)$-close to a solution $\phi$ to the perturbation-free dynamics $\mathcal{H}$. 
\end{thm}

{\it Proof}: According to Proposition \ref{Proposition dynamical properties} in Subsection 3.2, $\mathcal{H}$ satisfies the hybrid basic conditions, and is pre-forward complete from the compact set $\mathcal{C}$ since every $\phi \in \mathcal{S_H(C)}$ is complete (see Definition \ref{pre-forward completeness}). So from Lemma \ref{Proposition 6.34}, for every $\varepsilon>0$, $\tau \geq 0$, and $\rho : \mathbb{R}^N \rightarrow \mathbb{R}_{\geq 0}$, there exists a scalar $\sigma>0$ with the following property: for every solution $\phi_\sigma$ to $\mathcal{H}_{\sigma\rho}$ from $\mathcal{C}$, there exists a solution $\phi$ to $\mathcal{H}$ from $\mathcal{C}$ such that $\phi_\sigma$ and $\phi$ are $(\tau,\varepsilon)$-close, where $\mathcal{H}_{\sigma\rho}=(\mathcal{C},f_{\sigma\rho},\mathcal{D},G)$ is the $\sigma\rho$-perturbation of $\mathcal{H}$ and $f_{\sigma\rho}(x)=f(x)+\sigma\rho(x)\mathbb{B}=\omega \mathbf{1}_N+\sigma\rho(x)\mathbb{B}$ for every $x\in \mathcal{C}$. Note that if $p \in \sigma \rho(x) \mathbb{B}$, every solution $\phi_p$ to $\mathcal{H}_p$ from $\mathcal{C}$ is in fact the solution to $\mathcal{H}_{\sigma\rho}$, which implies that $\phi_p$ and $\phi$ are $(\tau,\varepsilon)$-close. \hfill$\blacksquare$

According to Theorem 2, the behavior of perturbed PCOs is close to the perturbation-free case, i.e., the solutions to the perturbed PCOs converge to the neighborhood of the synchronization set $\mathcal{A}$. Therefore, the phases of oscillators will remain close to each other under small frequency perturbations.

\section{Numerical Experiments}

\subsection{Unperturbed Case}
We first considered the unperturbed case, i.e., all oscillators had an identical frequency $\omega=2\pi$. 

First we considered $N=6$ PCOs on an undirected chain graph. Oscillators $1,\ldots,6$ adopted the PRFs (a), (b), (c), (d), (a), and (b) in Fig. \ref{fig:2}, respectively. The respective analytical expressions of these PRFs are given below.

\begin{equation}\label{PRF_a}
(a):  F_j(x_j)=\left\lbrace \begin{aligned}
& -0.6 x_j, & & \text{if} \ x_j \in [0, \, \pi) \\
& \big\{ -0.6\pi, \, 0.6\pi \big\}, & & \text{if} \ x_j=\pi \\
& \, 0.6(2\pi-x_j), & & \text{if} \ x_j \in (\pi, \, 2\pi] \\
\end{aligned}\right.
\end{equation}
\begin{equation}\label{PRF_b}
(b):  F_j(x_j)=\left\lbrace \begin{aligned}
& -0.7 x_j, & & \text{if} \ x_j \in [0, \, \frac{\pi}{2}) \\
& -0.35 \pi, & & \text{if} \ x_j \in [\frac{\pi}{2}, \, \pi) \\
& \big\{ -0.35 \pi, \, 0.35 \pi \big\}, & & \text{if} \ x_j=\pi \\
& \, 0.35 \pi, & & \text{if} \ x_j \in (\pi, \, \frac{3\pi}{2}] \\
& \, 0.7 (2\pi-x_j), & & \text{if} \ x_j \in (\frac{3\pi}{2}, \, 2\pi] \\
\end{aligned}\right.
\end{equation}
\begin{equation}\label{PRF_c}
(c): F_j(x_j)=\left\lbrace \begin{aligned}
&  -1.5\sin(0.5x_j), & & \text{if} \ x_j \in [0, \, \pi) \\
& \big\{ -1.5, \, 1.5 \big\}, & & \text{if} \ x_j=\pi \\
& \, 1.5\sin(0.5x_j), & & \text{if} \ x_j \in (\pi, \, 2\pi] \\
\end{aligned}\right.
\end{equation}
\begin{equation}\label{PRF_d}
(d): F_j(x_j)=\left\lbrace \begin{aligned}
&  -x_j^3/\pi^2+x_j^2/\pi-0.75x_j, \\
& \qquad \qquad \qquad \quad \ \text{if} \ x_j \in [0, \, \pi) \\
& \big\{ -0.75\pi, \, 0.75\pi \big\}, \\
& \qquad \qquad \qquad \quad \ \text{if} \ x_j=\pi \\
& \, -x_j^3/\pi^2+5x_j^2/\pi-8.75x_j+5.5\pi, \\
& \qquad \qquad \qquad \quad \ \text{if} \ x_j \in (\pi, \, 2\pi] \\
\end{aligned}\right.
\end{equation}

\begin{figure}[!h]
	\begin{center}
	\includegraphics[width=0.51\textwidth]{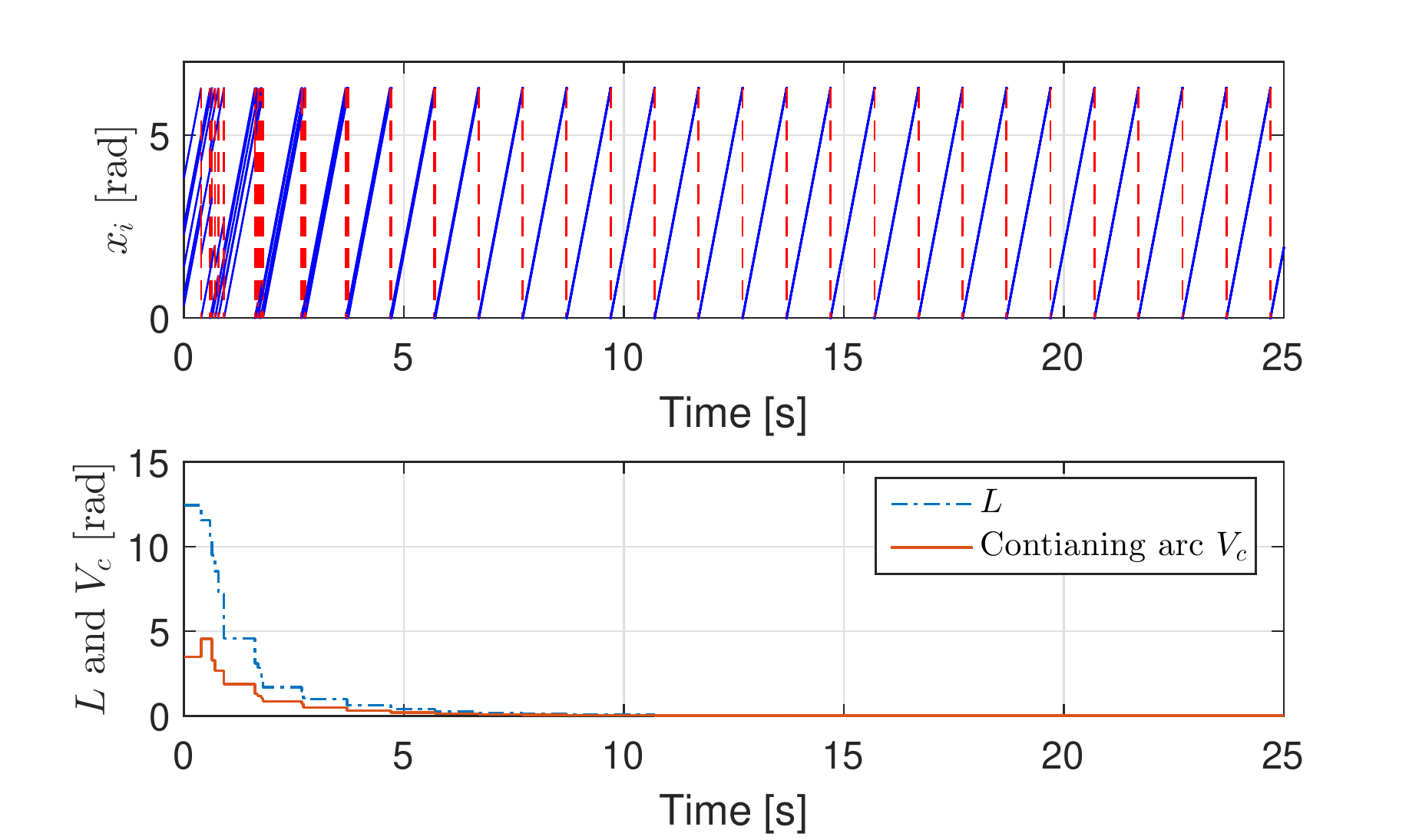}
	\end{center}
\caption{Evolutions of phases and $L$ for PCOs on an undirected chain graph.}
	\label{fig:6}
\end{figure}

The coupling strength $l_1,\ldots,l_6$ were set to $0.4$, $0.5$, $0.6$, $0.6$, $0.5$, and $0.4$, respectively. The initial phase $x(0,0)$ was randomly chosen from $\mathcal{C}\cup \mathcal{D}$. Fig. \ref{fig:6} shows the evolutions of phases and $L$. It can be seen that $L$ converged to $0$, which confirmed Theorem 1.

\begin{figure}[!h]
	\begin{center}
	\includegraphics[width=0.47\textwidth]{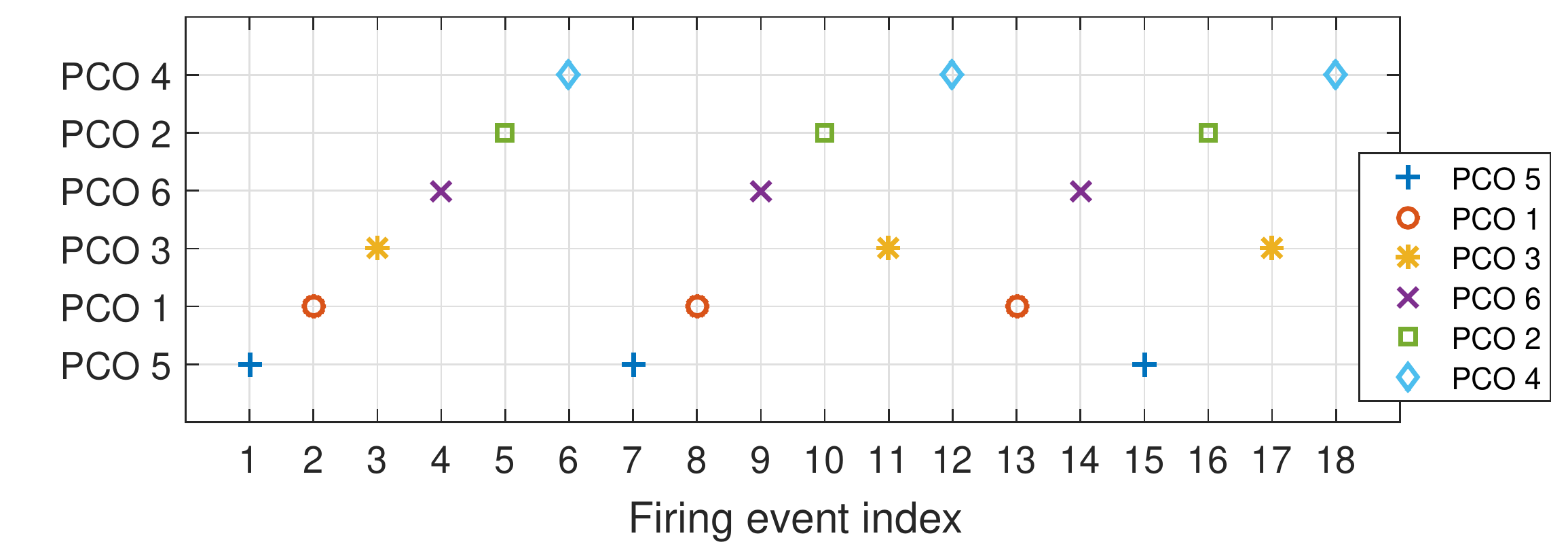}
	\end{center}
 \caption{Firing order of PCOs on the undirected chain graph.}
	\label{fig:07}
\end{figure} 
\begin{figure}[!h]
	\begin{center}
	\includegraphics[width=0.51\textwidth]{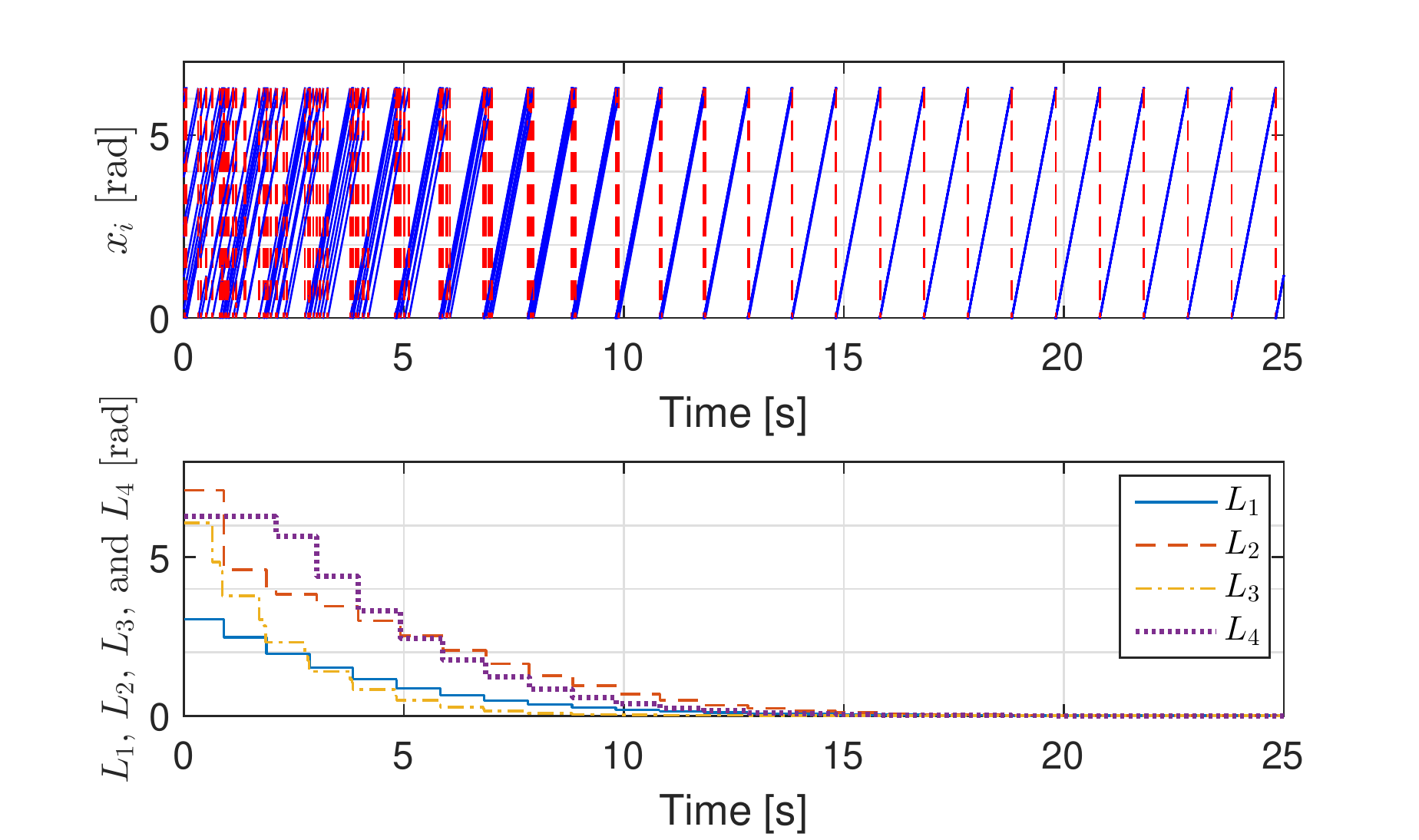}
	\end{center}
	\caption{Evolutions of phases, $L_1$, $L_2$, $L_3$, and $L_4$ for PCOs on a directed tree graph. PCOs synchronized as $L_1$, $L_2$, $L_3$, and $L_4$ converged to $0$.}
	\label{fig:10}
\end{figure} 

From the lower plot of Fig. \ref{fig:6}, we can also see that the length of the shortest containing arc $V_c$, which is widely used as a Lyapunov function in local synchronization analysis \cite{proskurnikov2016synchronization, kannapan2016synchronization, nunez2015synchronization, nunez2016synchronization}, is not appropriate for global PCO synchronization as it may not decrease monotonically. Along the same line, the firing order which is invariant in \cite{hong:05, goel2002synchrony, canavier2017globally}, and \cite{nunez2015synchronization}, is not constant in the considered dynamics as exemplified in Fig. \ref{fig:07}. These unique properties of chain and directed tree PCOs corroborate the novelty and importance of our results.

Then we considered $N=10$ PCOs on a directed tree graph, as illustrated in Fig. \ref{fig:1} (c). There are $4$ directed chains in this graph, namely, oscillators $1 \rightarrow 2 \rightarrow 5$, oscillators $1 \rightarrow 2 \rightarrow 4 \rightarrow 8$, oscillators $1 \rightarrow 3 \rightarrow 6 \rightarrow 9$, and oscillators $1 \rightarrow 2 \rightarrow 4 \rightarrow 7 \rightarrow 10$. The same as (\ref{Lyapunov}), $L_1$, $L_2$, $L_3$, and $L_4$ were defined to measure the degree of synchronization corresponding to the $4$ directed chains, respectively. Oscillators $1,\ldots,10$ adopted the PRFs (a), (b), (c), (d), (a), (b), (c), (d), (a), and (b) in Fig. \ref{fig:2}, respectively. The coupling strength $l_1,\ldots,l_{10}$ were set to $0.6$, $0.5$, $0.4$, $0.6$, $0.5$, $0.4$, $0.6$, $0.5$, $0.4$, and $0.6$, respectively. The initial phase $x(0,0)$ was randomly chosen from $\mathcal{C}\cup \mathcal{D}$. The convergence of $L_i$ ($i=1,\ldots,4$) to zero in Fig. \ref{fig:10} implies the synchronization of the $i$th directed chain, which confirmed Corollary 1. The simultaneous synchronization of all four directed chains also means synchronization of the entire directed tree graph, which confirmed Corollary 2.

\subsection{Perturbed Case}
\begin{figure}
	\begin{center}
	 \includegraphics[width=0.51\textwidth]{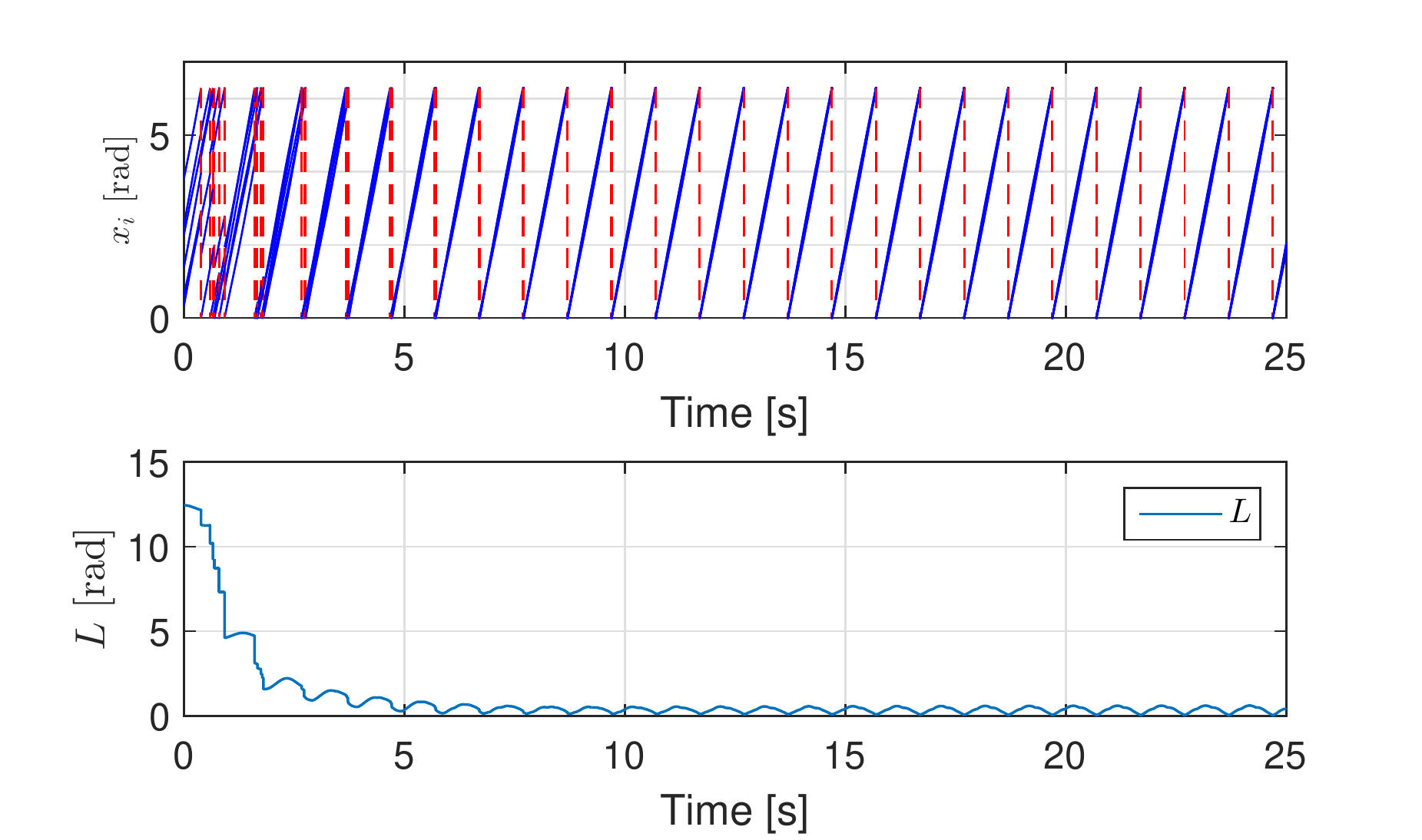}
	\end{center}
\caption{Evolutions of phases and $L$ for PCOs on an undirected chain graph under frequency perturbations.}
	\label{fig:11}
\end{figure}

We considered $N=6$ PCOs on an undirected chain graph with frequency perturbations on oscillator $k$ set to $p_k=0.5\sin(2\pi t+ 2\pi k/N)$. The other settings were the same as the undirected chain case. The evolutions of phases and $L$ were shown in Fig. \ref{fig:11}. It can be seen that the perturbed behaviors did not differ too much from the unperturbed case in Fig. \ref{fig:6}, and the solution converged to a neighborhood of the synchronization set $\mathcal{A}$ as $L$ approached a ball containing zero, which confirmed Theorem 2.

\section{Conclusions}

The global synchronization of PCOs interacting on chain and directed tree graphs was addressed. It was proven that PCOs can be synchronized from an arbitrary initial phase distribution under heterogeneous phase response functions and coupling strengths. The results are also applicable when oscillators are heterogeneous and subject to time-varying perturbations on their natural frequencies. Note that different from existing global synchronization results, the coupling strengths in our results can be freely chosen between zero and one, which is desirable since a very strong coupling strength, although can bring fast convergence, has been shown to be detrimental to the robustness of synchronization to disturbances. Given that a very weak coupling may not be desirable either due to low convergence speed which may allow disturbances to accumulate, the results give flexibility in meeting versatile requirements in practical PCO applications.


\begin{ack}
The authors would like to thank Francesco Ferrante for discussions and feedback which greatly strengthened the paper.
\end{ack}


\section*{Appendix A:  Lemma \ref{Lemma2}}

\begin{Lemma 2}\label{Lemma2}
	For $N$ PCOs interacting on an undirected chain, if the PRF $F_j(x_j)$ satisfies Assumption 3 and $l_j\in (0, \, 1)$ holds for all $j \in \mathcal{V}$, then $L$ in (\ref{Lyapunov}) cannot be retained at any nonzero value along a complete solution $\phi$.
\end{Lemma 2}

{\it Proof}: We use proof of contradiction. Since $L\in [0, N\pi]$ holds, we suppose that for some $r \in (0, \, N\pi]$, $L$ is retained at $r$ along a complete solution $\phi$. From Lemma \ref{Lemma1}, to keep $L$ at $r$, we must have
\begin{equation}\label{Lemma_21}
\Delta_{i-1}^+ =\Delta_{i-1}-\delta_{i-1}, \quad \Delta_{i-2}^+ =\Delta_{i-2}+\delta_{i-1}
\end{equation}
or
\begin{equation}\label{Lemma_22}
\Delta_{i}^+ =\Delta_{i}-\delta_{i+1}, \quad \Delta_{i+1}^+ =\Delta_{i+1}+\delta_{i+1}
\end{equation}
if the left-neighbor oscillator $i-1$ or right-neighbor oscillator $i+1$ exists when oscillator $i$ fires, respectively. Next we show that $\Delta_{N}$ will exceed $\pi$, which contradicts the constraint $0\leq \Delta_{i}\leq \pi$ for $i\in \mathcal{V}$.

Given $1 \notin \mathcal{N}_i^{out}$ and $N \notin \mathcal{N}_i^{out}$ for $i=3,4,\ldots,N-2$, both $x_1^+=x_1$ and $x_N^+=x_N$ hold when oscillators $3,4,\ldots,N-2$ fire, which leads to $\Delta_N^+=\Delta_N$. Similarly, $N \notin \mathcal{N}_1^{out}$ (resp. $1 \notin \mathcal{N}_N^{out}$) implies $x_N^+=x_N$ (resp. $x_1^+=x_1$) when oscillator $1$ (resp. $N$) fires, which leads to $\Delta_N^+=\Delta_N$ when oscillator $1$ or $N$ fires.

So we focus on the evolution of $\Delta_N$ when oscillators $2$ and $N-1$ fire. According to Lemma \ref{Lemma3} in Appendix B, neither oscillator $2$ nor oscillator $N-1$ will stop firing. Without loss of generality, we assume that oscillator $2$ fires at time $(t_2^*, \, k_2^*)$. From (\ref{Lemma_21}) we have $\Delta_{N}^+ =\Delta_{N}+\delta_{1}$. Similarly, from (\ref{Lemma_22}) we have $\Delta_{N}^+ =\Delta_{N}+\delta_{N}$ when oscillator $N-1$ fires. Since $\delta_{1}$ and $\delta_{N}$ are nonnegative, we have $\Delta_{N}^+ \geq \Delta_{N}$. To prove that $\Delta_N$ will surpass $\pi$, we need to show that at least one of the following statements is true: 
\begin{enumerate}[1)]
	\item $\delta_{1} = 0$ cannot always hold when oscillator $2$ fires;
	\item $\delta_{N} = 0$ cannot always hold when oscillator $N-1$ fires.
\end{enumerate}

\begin{figure}
	\begin{center}
		\includegraphics[width=0.25\textwidth]{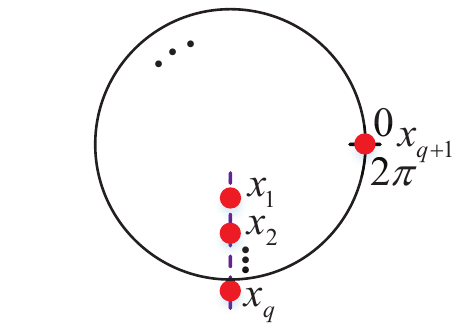}
	\end{center}
	\caption{Illustration of a set of $q \geq 2$ neighboring oscillators being synchronized.}
	\label{fig:6_2}
\end{figure}

Proof of statement 1): Given $l_1\in (0, \, 1)$, according to (\ref{PRF_detail}) and (\ref{ddelta_i-1}), $\delta_{1} = 0$ holds if and only if $x_1 = 0$ or $x_1 = 2\pi$ holds, which means that oscillators $1$ and $2$ are synchronized when oscillator $2$ fires. So we need to show that oscillators $1$ and $2$ cannot always be synchronized when oscillator $2$ fires. More generally, we assume that there is a set of $q \geq 2$ oscillators $1,2,\ldots,q$ being synchronized and having phases different from oscillator $q+1$. According to Lemma \ref{Lemma3} in Appendix B, oscillator $q+1$ will not stop firing in this situation. We assume that oscillator $q+1$ fires at time $(t^*_{q+1}, \, k^*_{q+1})$, and $x_{1}=\ldots=x_q \in [\pi, \, 2\pi)$ holds when oscillator $q+1$ fires, as illustrated in Fig. \ref{fig:6_2}. Note that the case of $x_{1}=\ldots=x_q \in (0, \, \pi]$ can be proved by following the same line of reasoning. Given $0<l_q<1$, from (\ref{PRF_detail}) and (\ref{ddelta_i-1}) we have $0<\delta_q<2\pi-x_q$. Since oscillator $q$ is the left-neighbor of oscillator $q+1$, according to (\ref{Lemma_21}), when oscillator $q+1$ fires we have $\Delta_{q}^+ =\Delta_{q}-\delta_{q} =2\pi-x_q-\delta_{q} >0$ and $\Delta_{q-1}^+ =\Delta_{q-1}+\delta_{q}=0+\delta_{q}>0$. So oscillator $q$ escapes from the set of synchronized oscillators due to $\Delta_{q-1}^+>0$ and will fire next. Similarly, when oscillator $q$ fires, the left-neighbor oscillator $q-1$ will escape from the set of synchronized oscillators and fire next. Iterating this argument, when oscillator $3$ fires, the left-neighbor oscillator $2$ will escape from the set of synchronized oscillators and fire next. So we have $x_{2} \neq x_1$, i.e., oscillators $1$ and $2$ are not synchronized when oscillator $2$ fires. Therefore, $\delta_1 = 0$ cannot always hold when oscillator $2$ fires. 

Similarly, we can prove statement 2), i.e., $\delta_N$ cannot always be $0$ when oscillator $N-1$ fires, and thus $\Delta_N$ will keep increasing. Since $\delta_1$ and $\delta_N$ will not converge to $0$ unless synchronization is achieved, $\Delta_N$ will surpass $\pi$, which contradicts the constraint $0\leq \Delta_{i}\leq \pi$ for $i\in \mathcal{V}$. Therefore, $L$ cannot be retained at any nonzero value along a complete solution $\phi$. \hfill$\blacksquare$

\section*{Appendix B:  Lemma \ref{Lemma3}}

\begin{Lemma 2}\label{Lemma3}
	For $N$ PCOs interacting on an undirected chain, if the PRF $F_j(x_j)$ satisfies Assumption 3 and $l_j\in (0, \, 1)$ holds for all $j \in \mathcal{V}$, we have the following results:
	\begin{enumerate}[1)]
		\item Neither oscillator $2$ nor oscillator $N-1$ will stop firing;
		\item Oscillator $q+1$ will not stop firing if oscillators $1,\ldots,q$ ($2\leq q \leq N-1$) have been synchronized and oscillator $q+1$ is not synchronized with these $q$ oscillators. Similarly, oscillator $N-q$ will not stop firing if oscillators $N-q+1,\ldots,N$ have been synchronized and oscillator $N-q$ is not synchronized with these $q$ oscillators.
	\end{enumerate}
\end{Lemma 2}

{\it Proof}: We first use proof of contradiction to prove statement 1). Suppose that oscillator $2$ stops firing after time instant $(t'_{2}, \, k'_2)$, then $x_{2}$ will stay in $[0,\pi]$. This is because if $x_{2} \in (\pi, 2\pi)$ holds, it will evolve continuously to $2\pi$ and fire, and receiving pulses from other oscillators can only expedite this process under the PRFs in Assumption 3. Since oscillator $2$ only receives pulses from oscillators $1$ and $3$, without loss of generality, we suppose at time $(t'_1, \, k'_1)$ that oscillator $1$ fires and resets its phase to $0$. Note that oscillator $1$ will fire at a period of $T_1=2\pi/\omega$ since its only neighbor oscillator $2$ stops firing. After receiving the pulse, oscillator $2$ updates its phase to $x_2^+=x_2+l_2 F_2^{(1)}(x_2) \in [0, \pi)$. If oscillator $2$ does not receive any other pulse before its phase surpasses $\pi$, it will fire, which contradicts the assumption. So we suppose that oscillator $3$ fires at time $(t'_3, \, k'_3)$ before $x_2$ surpasses $\pi$, which implies $t'_3-t'_1 \leq \pi/\omega$. Since the time it takes for phase evolving from $0$ to $\pi$ is at least $\pi/\omega$ and after reaching $\pi$ oscillator $3$ will not fire immediately even if it receives a pulse under given PRFs and coupling strengths, the length of oscillator $3$'s firing period $T_3$ satisfies $T_3 > \pi/\omega$. There are two cases in this situation, $t'_1=t'_3$ and $t'_1<t'_3$, respectively:

Case 1: If $t'_1=t'_3$ holds, then the length of time interval for oscillator $2$ receiving the next pulse after $(t'_3, \, k'_3+1)$ is greater than $\pi/\omega$. Since $x_2(t'_3, \, k'_3+1)\geq 0$ holds, $x_2$ will be greater than $\pi$ when receiving the next pulse. So oscillator $2$ will fire again, which contradicts the assumption.

Case 2: If $t'_1<t'_3$ holds, then we have $x_2(t'_3, \, k'_3+1) > 0$ due to $x_2(t'_3, \, k'_3)=x_2(t'_1, \, k'_1+1)+\omega(t'_3-t'_1) > 0$ under given PRFs and coupling strengths. Since $t'_3-t'_1 \leq \pi/\omega$ holds, after time interval $[\pi-x_2(t'_3, \, k'_3+1)]/\omega$ which is less than $\pi/\omega$, we have $x_1<2\pi$, $x_3<\pi$, and $x_2=\pi$. So $x_2$ will be greater than $\pi$ when receiving the next pulse, and thus oscillator $2$ will fire again, which contradicts the assumption.

Therefore, oscillator $2$ will not stop firing. Similarly, we can prove that oscillator $N-1$ will not stop firing either. 

Next we prove statement 2). Suppose that oscillator $q+1$ stops firing after time $(t'_{q+1}, \, k'_{q+1})$. Since oscillators $1,\ldots,q$ will not receive any pulses from other oscillators, they will remain synchronized and oscillator $q$ will fire with a period of $T_q=2\pi/\omega$. The same as statement 1), the length of oscillator $q+2$'s firing period $T_{q+2}$ satisfies $T_{q+2} > \pi/\omega$ and oscillator $q+1$ will not stop firing if oscillator $q+1$ has a phase different from synchronized oscillators $1,\ldots,q$. Similarly, we can prove that oscillator $N-q$ will not stop firing either if oscillator $N-q$ has a phase different from synchronized oscillators $N-q+1,\ldots,N$. \hfill$\blacksquare$

\bibliographystyle{unsrt}
\bibliography{abbr_bibli}

\end{document}